\documentclass{tip}
\usepackage{orcidlink}
\newtheorem{definition}{Definition} 
\newtheorem{theorem}[definition]{Theorem}
\newtheorem{lemma}[definition]{Lemma}
\newcommand{\autora}{Rico~Zöllner}
\newcolumntype{K}[1]{>{\centering\arraybackslash}p{#1}}
\newcommand{\mc}[1]{\multicolumn{2}{K{1.8cm}}{#1}}
\allowdisplaybreaks
\usepackage{tabulary}

\title{\vspace{-2.5cm}\hrulefill\\\color{darkblue}\sffamily\LARGE\bfseries A New Perspective on Double-S Curve Motions of Higher Order and Optimal Motion Planning\\}
\author{\sffamily\bfseries \autora\orcidlink{0000-0002-3544-6622}${}^a$}
\affil{\sffamily${}^a$Institute of Material Handling and Industrial Engineering, TU~Dresden,\newline01062 Dresden, Germany}
\date{\sffamily\today\\\hrulefill}

\usepackage{minted}
\usepackage{listings}

\definecolor{vscode-bg}{RGB}{30,30,30}
\definecolor{vscode-text}{RGB}{220,220,220}
\definecolor{vscode-keyword}{RGB}{86,156,214}
\definecolor{vscode-string}{RGB}{206,145,120}
\definecolor{vscode-comment}{RGB}{106,153,85}
\definecolor{vscode-func}{RGB}{220,220,170}

\lstdefinestyle{vscode}{
	backgroundcolor=\color{vscode-bg},
	basicstyle=\ttfamily\small\color{vscode-text},
	keywordstyle=\color{vscode-keyword}\bfseries,
	stringstyle=\color{vscode-string},
	commentstyle=\color{vscode-comment},
	identifierstyle=\color{vscode-text},
	emphstyle=\color{vscode-func},
	numbers=left,
	numberstyle=\tiny\color{gray},
	stepnumber=1,
	numbersep=10pt,
	showstringspaces=false,
	breaklines=true,
	tabsize=4
}
\begin{document}
	\maketitle
	\thispagestyle{plain}
	\vspace{-2.0em}
	
	\parbox[t]{0.30\linewidth}{\textbf{Keywords}
		\begin{flushleft}
			\begin{itemize}
			\item Optimal Motion Planning
			\item Double-S Curve Trajectory
			\item Minimum-Time Motion
			\item High-Order Control
				\end{itemize}
			\end{flushleft}}
	\hfill
	\parbox[t]{0.60\linewidth}{\textbf{Abstract}\\{\footnotesize This paper presents and proves an equation for the time horizon of symmetric trajectories with zero boundary conditions and bounded derivatives of arbitrary order. This equation holds regardless of the number of phases comprising the associated motion, which avoids case distinctions in calculations. Application examples of motions with minimum time, minimum velocity, and minimum acceleration are discussed. Furthermore, an algorithm is derived that reduces the time minimization problem to solving a system of equations. This algorithm avoids nested case distinctions and complex optimizations.} }
	
	\vspace{2.0em} 
	\noindent\rule{\textwidth}{0.4pt}
	\vspace{-3.0em}
	
	\section{Introduction}\label{Sec1}
Optimal motion planning has been done for many years, and it is still an active research field -- particularly from the mathematical-algorithmical and the numerical-computational perspective. Some important applications are robotics (see~\cite{Zhou,Mohanan}), aeronautics and astronautics (see~\cite{Trelat,Izzo}), or transport planning (see~\cite{Katrakazas}) -- just to mention a few. A rather simple but important task is to find the minimum-time motion for accomplishing a given distance~$s$ under certain restrictions: e.g., the velocity is bounded, or the acceleration is bounded, or both (see~\cite{Biagiotti} for an excellent book about the kinematics of trajectory planning); surveys in the field of optimization-based motion planning are, e.g.,~\cite{Almasri,yang19}. Common approaches consider the jerk as third derivative as well~\cite{Alpers,Lee,R1}. Such approaches work in a 7-phase motion (SPM) profile (which is also referred as 7-phase motion or double-S curve, see~\cite{Castain} for an early work and \cite{Lu16,Mu09} as well) to find relations between the minimal time, the distance and the boundary values for several cases and under certain boundary conditions. Fourth order motions, i.e. snap-controlled motions with constraints for the first four derivations of the position, are considered in~\cite{Lambrechts,Kai}. Particularly, double-S (or simply S-curve) motion planning remains a vital topic in high-precision robotic trajectory generation: for instance \cite{liu23} propose a time-optimal asymmetric S-curve trajectory for redundant manipulators under kinematic constraints, using a whale-optimization algorithm to schedule jerk-limited motion phases while respecting both end-effector and joint limits. In parallel, \cite{li23} introduce a general asymmetric S-curve profile optimized to minimize execution time while keeping residual vibration low, demonstrating analytic performance gains over symmetric or simpler S-curve profiles. Moreover, \cite{liu23} present a time-optimal, jerk-continuous trajectory planning method for robot manipulators that enforces third-order kinematic constraints (velocity, acceleration, jerk) to produce smooth, practically trackable motions. In industrial motion-planning scenarios, \cite{lin18} propose an SPM optimizer under joint kinematics limits and waypoint constraints, framing the trajectory planning as a nonlinear optimization to minimize cycle time. 
Recent work also pushes beyond purely time-optimal planning: \cite{bao24} develop a multi-objective optimal trajectory planner for manipulators that jointly optimizes time and jerk (or smoothness), using quintic B-splines and a meta-heuristic to respect kinematic constraints while trading off between speed and smoothness.

The aim of this paper is to reveal new perspectives on symmetric double-S curves of arbitrary order. For a treatment of asymmetric trajectories see~\cite{Ha}. Following the work of~\cite{Nguyen} (see also~\cite{Fang,Boryga}), a general equation for the duration of the motion $T$ is proven. In the case of jerk-controlled motion\footnote{The term ``jerk-controlled'' means that the velocity, the acceleration and the jerk are assumed to be bounded.}, this equation is
\begin{equation}
	T = \frac{s}{v} + \frac{v}{a} + \frac{a}{j}.
	\label{Eins}
\end{equation}
Here, $v$ denotes the maximum velocity achieved during the motion; $a$ represents the maximum acceleration achieved; and $j$ stands for the corresponding jerk. The generalization to higher order appears quite obvious. Furthermore, it is noteworthy that this Equation~(\ref{Eins}) holds regardless of how many phases the motion actually consists of\footnote{For example, a jerk-controlled motion can consist of four, five, six or seven phases. Considering the next order (snap), there can occur 15 phases -- see~\cite{Biagiotti} where an analysis is done phase by phase.}.

Moreover, since the equation is true for all symmetric\footnote{I.e. roughly speaking acceleration = deceleration, see below.} trajectories with vanishing boundary values, it can be used to describe not only the time-minimal motion over the distance~$s$. Conversely, if a time horizon and a distance are given, the question arises as to which trajectory -- among all possible motions -- minimizes, for example, the achieved velocity (or acceleration, or jerk, etc.). Equation~(\ref{Eins}) and its generalization are also suitable for such questions. In the following, a systematic approach to various optimizations and arbitrary orders is developed using a single equation.

The paper is organized as follows. Section~\ref{Section2} precisely defines terms such as ``symmetric'', which were only used intuitively in the introduction. The class of admissible trajectories is defined and the optimization tasks under consideration are set up. To familiarize ourselves with the topic, Section~\ref{Section3} deals with the jerk-controlled motion, and Equation~(\ref{Eins}) is proven on an elementary level for the four possible cases by case distinctions. This is the common approach, done e.g., in~\cite{Biagiotti,Alpers}. Section~\ref{Section4} presents the generalization of Equation~(\ref{Eins}) and the general proof for arbitrary orders. Section~\ref{Section5} is devoted to applications with the focus lying on optimal motion planning, especially for jerk- and snap-controlled movements. In addition, the solution method for finding the minimum-time motions is explained. A brief summary is provided in Section~\ref{Section6}.
	\section{Problem Setup} \label{Section2}
Let us consider one-dimensional trajectories encoded in $x(t)$ where $x \in \mathbb{R}$ denotes the position as a function of the time $t$ which ranges in the interval $[0,T]$ with $T > 0$. In addition, the order parameter is denoted by $N$ with $N \in \mathbb{N}$ and $N > 0$. We assume $x(t)$ to be $(N-1)$-times continuously differentiable w.r.t. $t$, and the $N$-th derivative of $x(t)$ should be piecewise continuous with countable many finite jumps. Then, $\dot{x}(t)$ is the velocity, and -- if existing -- $\ddot{x}(t)$ and $\dddot{x}(t)$ stand for the acceleration and the jerk, respectively\footnote{According to Newton, the dot indicates derivatives w.r.t. $t$. And in general, the $k$-th derivative w.r.t. $t$ is written as $x^{(k)}(t).$}.

\begin{definition}
	A trajectory $x(t)$ of order $N$ and with time horizon $T$ is defined to be \textup{admissible}, if it fullfils the following boundary conditions and constraints:
\begin{enumerate}[label=(\roman*)]
	\item $x^{(n)}(0) = 0$ for all $0\leq n \leq N-1$.
	\item $x^{(n)}(T) = 0$ for all $1\leq n \leq N-1$.
	\item $x(T) = s$ with a distance parameter $s>0$.
	\item For all $t\in [0,T]$ and $1 \leq n \leq N$ one has $\left|x^{(n)} (t) \right| \leq w_n$ with bounds $w_n > 0$. \label{Def1iv}
\end{enumerate} \label{Def1}
\end{definition}

In what follows, we have to distinguish between the given bounds $w_n$ which are sometimes more concretely denoted by $v_{\max}, a_{\max}, j_{\max}$ etc. and actually achieved values  $\max\limits_{t\in[0,T]} \left|x^{(n)} (t) \right|$. These values will be denoted by $v,a,j$ etc. for $n=1,2,3$, respectively. Some authors call the letters ones ``peak values'' (in contrast to ``maximum values'' for $w_n$) which may be quite misleading if -- for example -- the velocity $v$ is held for a certain time interval (cruising phase).

We are now prepared to formulate two external problems:\\

\noindent\textbf{Minimum-Time Trajectory (MTT)}

\noindent\textit{Let $N$ and $s$ be given as well as $(w_n)^N_{n=1}$. Find the minimal value of $T$ such that an admissible trajectory exists.}\\

Cleary, this is nothing else than the classical problem of finding the minimum-time motion for accomplishing a distance $s$ with bounded velocity $(N=1)$, acceleration $(N=2)$, jerk $(N = 3)$, $\ldots$\\

\noindent\textbf{Minimum-Derivative Trajectory (MDT)}

\noindent\textit{Let $s,~T$ and $N>2$ be given as well as $(w_n)^N_{n=1}$. For an order $M$ with $1 \leq M < N$, find an admissible trajectory such that
\begin{equation}
	\max\limits_{t\in[0,T]} \left|x^{(M)} (t) \right| \label{max}
\end{equation}
becomes minimal.}\\

Here, both distance $s$ and the time horizon $T$ are given\footnote{Clearly, the given time horizon of MDT has to be larger than the one of the minimum-time solution of MTT to ensure the existence of a solution.}, and the task is to find a trajectory, where a certain derivative (e.g., the velocity if $M=1$) is minimized. Such minimum-velocity motions or minimum acceleration motions (i.e. $M=2$ for $N \geq 3$) are considered because such trajectories protect the device since they avoid wear and tear or large forces for instance.

For both MTT and MDT the key characteristics of the optima are known~\cite{R1}. Basically, the $N$-th derivative of an optimal trajectory is piecewise constant\footnote{This is due to Pontryagin's maximum principle, see~\cite{pontryagin} for instance.}, and optimal trajectories enjoy a certain kind of symmetry in such a way that -- for example -- accelerating and decelerating are time-reversed processes. The following definition completes this section by formalizing the above mentioned relevant trajectory type.

\begin{definition} Let $x$ be an admissible trajectory of order $N$ and with time horizon $T$.
\begin{enumerate}[label=(\roman*)]
	\item An interval $I \subseteq [0,T]$ where $x^{(n)}(t) \neq 0$ for all  $t \in I$ is called \textup{supporting interval} of $x^{(n)}$ (with $1 \leq n \leq N$).
	\item A \textup{maximal supporting interval} $I$ of $x^{(n)}$ is a supporting interval $I$ of $x^{(n)}$ such that there is no supporting interval $J$ with $I$$ \subsetneq J \subseteq [0,T]$.
	\item The function $x^{(n)}$ (with $1 \leq n \leq N$) is called symmetric provided that for every maximal supporting interval $I = \left]a,b \right[$ one has
	\begin{equation}
		x^{(n)} (a + t) = x^{(n)}(b-t)
	\end{equation}
	for all $0 \leq t \leq b-a$. The trajectory $x$ is defined to be \textup{symmetric} if every order $x^{(n)}$ with $1 \leq n \leq N$ is symmetric. \label{Def2iii}
	\item A trajectory $x$ is defined to be \textup{pre-optimal} provided it is symmetric and there exists a constant $c > 0$ such that $x^{(N)}  (t) \in \left\lbrace -c,0,c\right\rbrace$ for all $t \in [0,T]$. \label{Def2iv}
\end{enumerate} \label{Def2}
\end{definition} 
A supporting interval is an interval, where the respective derivative does not vanish. If a supporting interval is maximal it can not be extended to the left or the right. That means that the boundary values of the respective derivative are zero for a maximal supporting interval. In Figure~\ref{Figure1}, right panel, any subinterval of $[0,T]$ is a supporting interval of $\dot{x}$ but only $[0,T]$ is a maximal supporting interval. Regarding the acceleration $\ddot{x}$, the maximal supporting intervals are $[0,t_1]$ and $[T-t_1,T]$. The proposed symmetry of Definition~\ref{Def2}~\ref{Def2iii} can be seen from that panel as well for $\dot{x}$. This symmetry means that any order restricted to a maximal supporting interval is axisymmetric to a vertical axis drawn through the middle of the interval. Definition~\ref{Def2}~\ref{Def2iv} means that the highest order is piecewise constant. Obviously, the constant $c$ in Definition~\ref{Def2}~\ref{Def2iv} should be less than or equal to the bound $w_N$ (see Definition~\ref{Def1}~\ref{Def1iv}).

The next section recalls briefly the pre-optimal trajectories up to the 7-phase motion at an elementary level.
\section{Seven-Phase Motions}\label{Section3}
Before we tackle the general proof, we briefly recall pre-optimal trajectories within the SPM framework while climbing from $N=1$ to $N=3$.

Note that the case of first-order motions (i.e. $N=1$) is trivial since the only relevant equation is 
\begin{equation}
	T = \frac{s}{v},
\end{equation}
where $v$ can either stand for the maximal velocity (i.e. $v = w_1$) in the case of MTT, or $v$ stands for the minimal possible velocity in case of MDT. In summary, we have a uniform motion in any case.
\subsection{Acceleration-Controlled Motions}\label{3.1}
\begin{figure}[h!]
	\centering
	\includegraphics[width=0.4\textwidth]{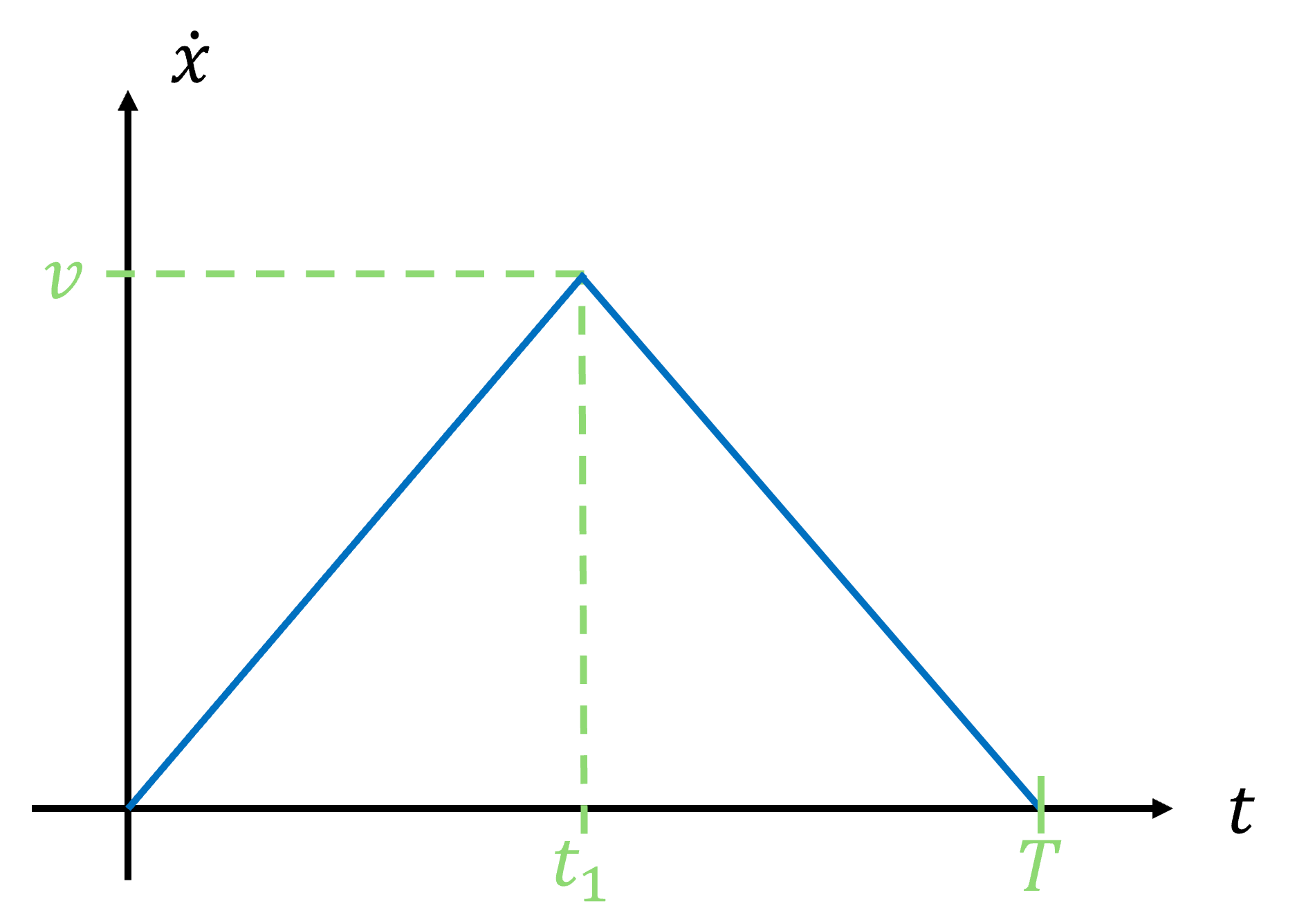}
	\includegraphics[width=0.55\textwidth]{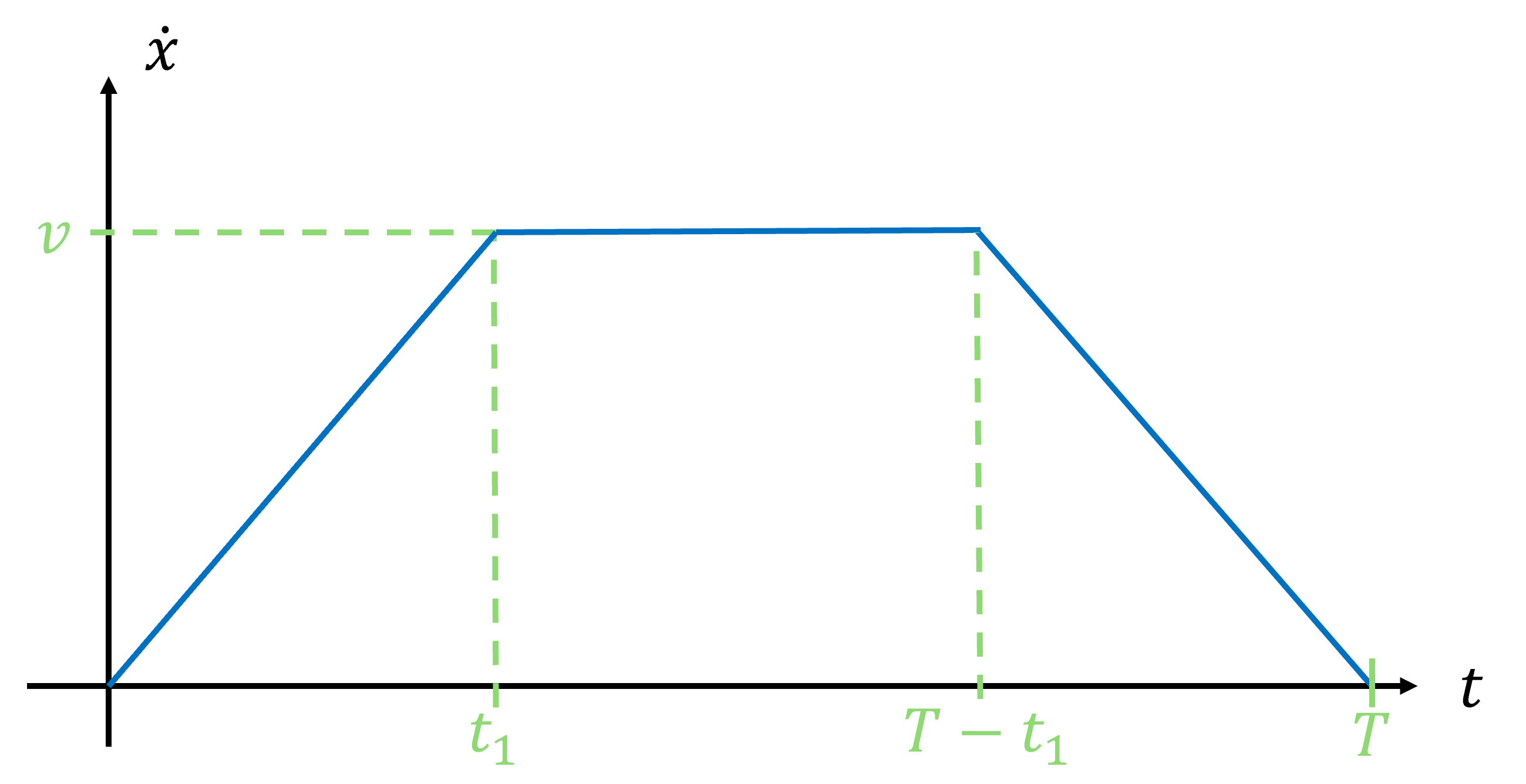}
	\caption{Velocity $\dot{x}$ as a function of the time $t$ for a triangular profile (left) and a trapezoidal profile (right). The maximal reached velocity is denoted by $v$, i.e.~$v~=~\max\limits_{t\in [0,T]} \left|\dot{x}(t) \right|$.} \label{Figure1}
\end{figure}

If we consider second-order motions ($N=2$), the acceleration enters the game. Now there are two cases of pre-optimal trajectories depending on the involved parameters $T$, $s$ and $w_1$, $w_2$: The velocity profile $\dot{x}(t)$ is triangular or trapezoidal (see Figure~\ref{Figure1} left and right, respectively). The constant $c$ of Definition~\ref{Def2}~\ref{Def2iv} is denoted with $a$ since it has the meaning of an acceleration, i.e. $\ddot{x} (t) \in \left\lbrace -a,0,a\right\rbrace$ for all $t \in [0,T]$, and $a \leq w_2$.

Let $v$ denote the maximum velocity reached during the motion\footnote{Note again, that achieved maximum of the velocity does not need to be equal to the given bound $w_1 = v_{\max}$; we always have $v \leq w_1$.}, i.e. $v = \max\limits_{t\in [0,T]} |\dot{x}(t)|$. In the triangular case we have
\begin{equation}
	v = a  t_1,~s = v  t_1 \Rightarrow t_1 = \frac{v}{a} = \frac{s}{v}
\end{equation}
and $T = 2t_1 = t_1 + t_1 = \frac{s}{v} + \frac{v}{a}$. Analogously, in the trapezoidal case we have again $v = a  t_1$ but $s = vt_1 + v  (T-2t_1) = v  (T-t_1)$ which again leads to 
\begin{equation}
	T = \frac{s}{v} + \frac{v}{a}. \label{T=}
\end{equation}
This is the tiny sister of Equation~(\ref{Eins}). In Section~\ref{Section5}, we turn to the solutions of MTT and MDT based on Equation~(\ref{T=}).

\subsection{Jerk-Controlled Motions}\label{3.2}
\begin{figure}[h!]
	\centering
	\includegraphics[width=0.4\textwidth]{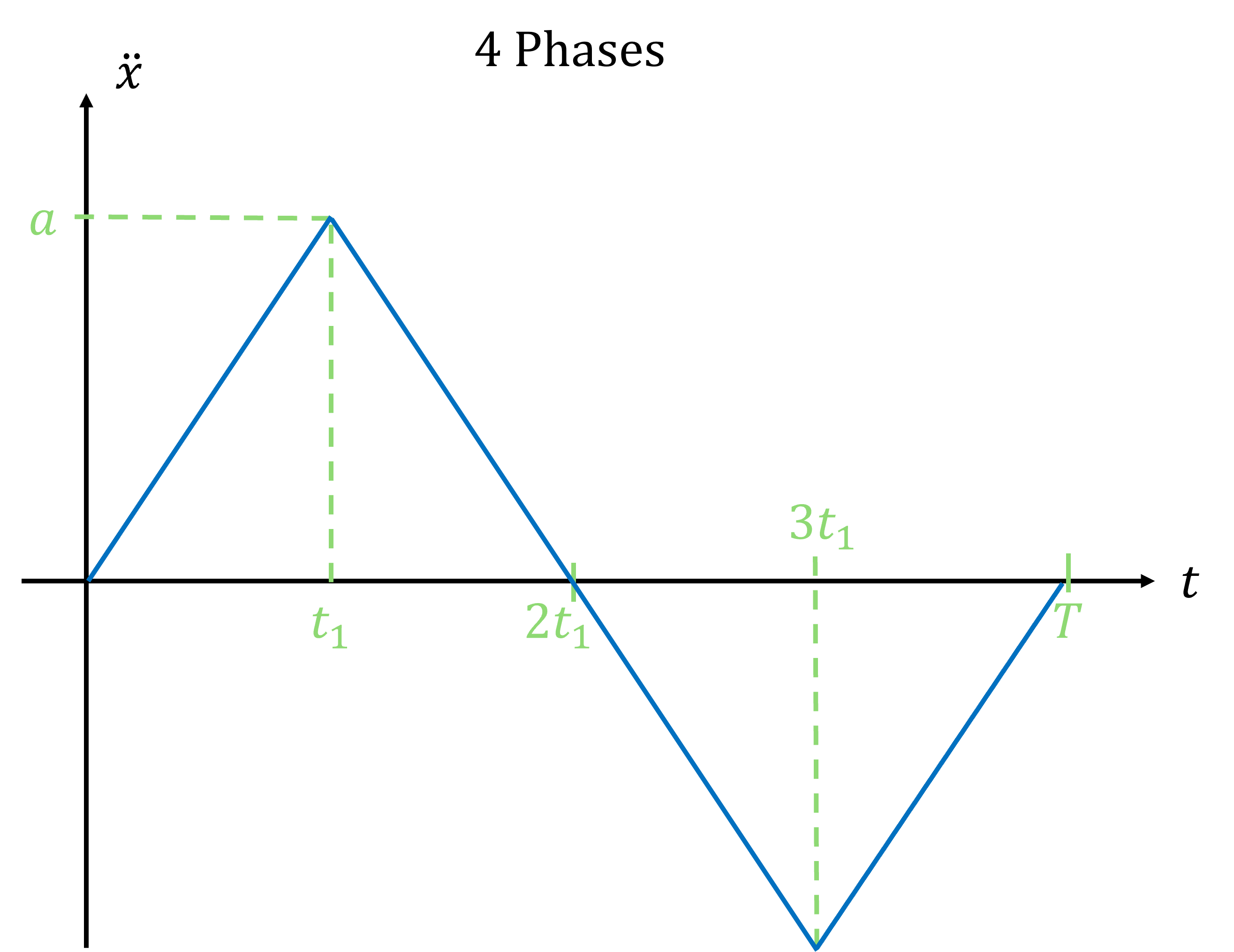}
	\includegraphics[width=0.4\textwidth]{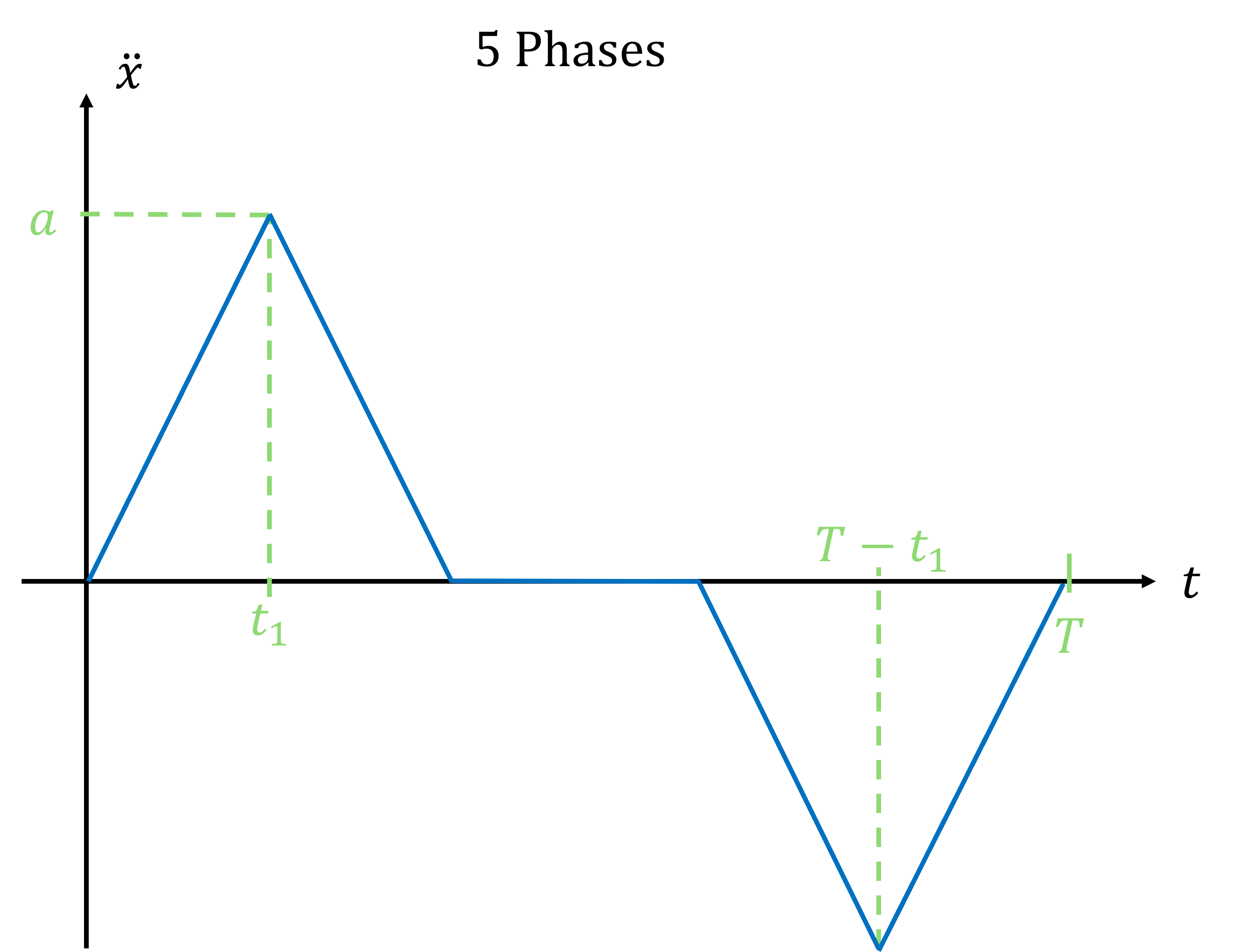}\\
	\includegraphics[width=0.4\textwidth]{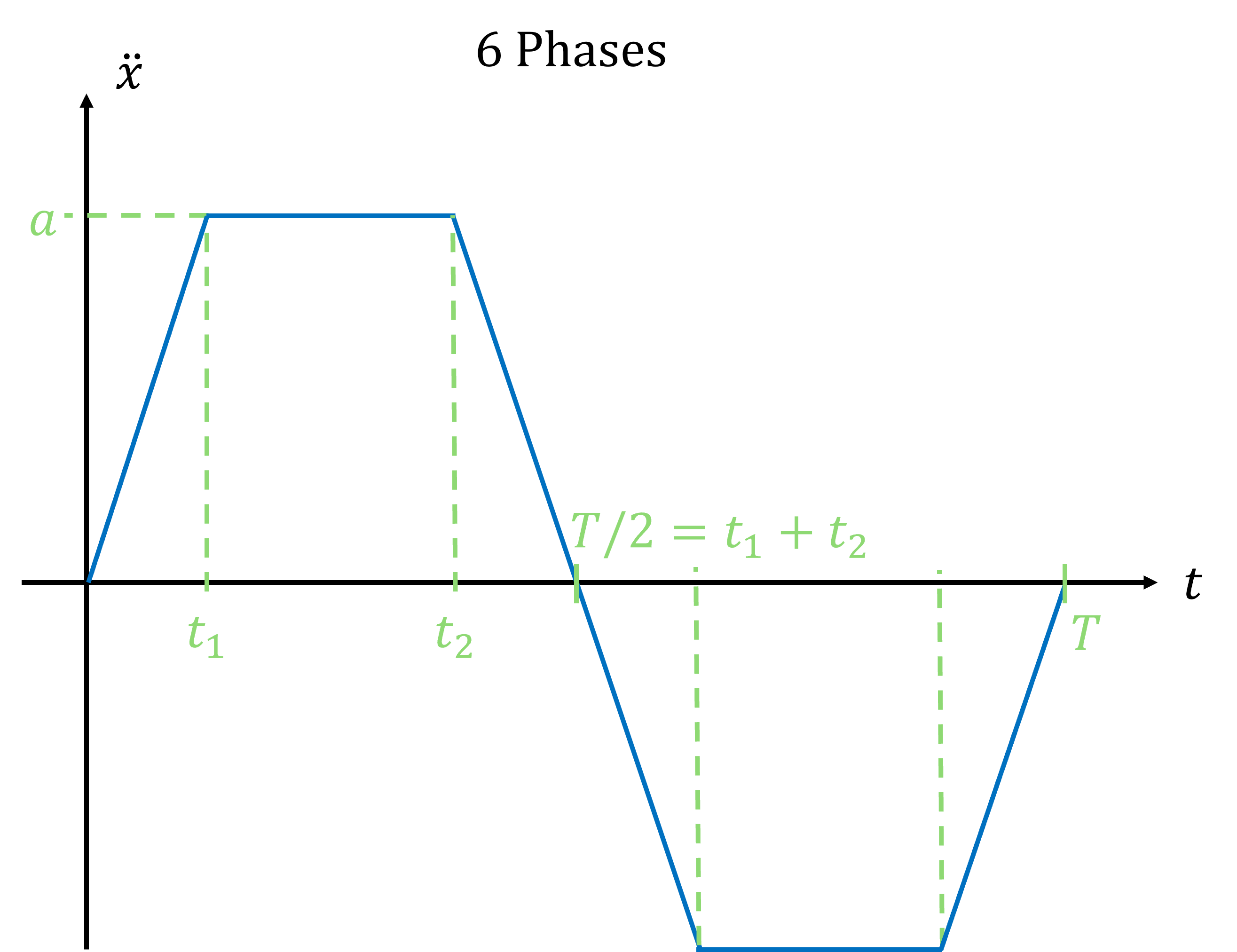}
	\includegraphics[width=0.4\textwidth]{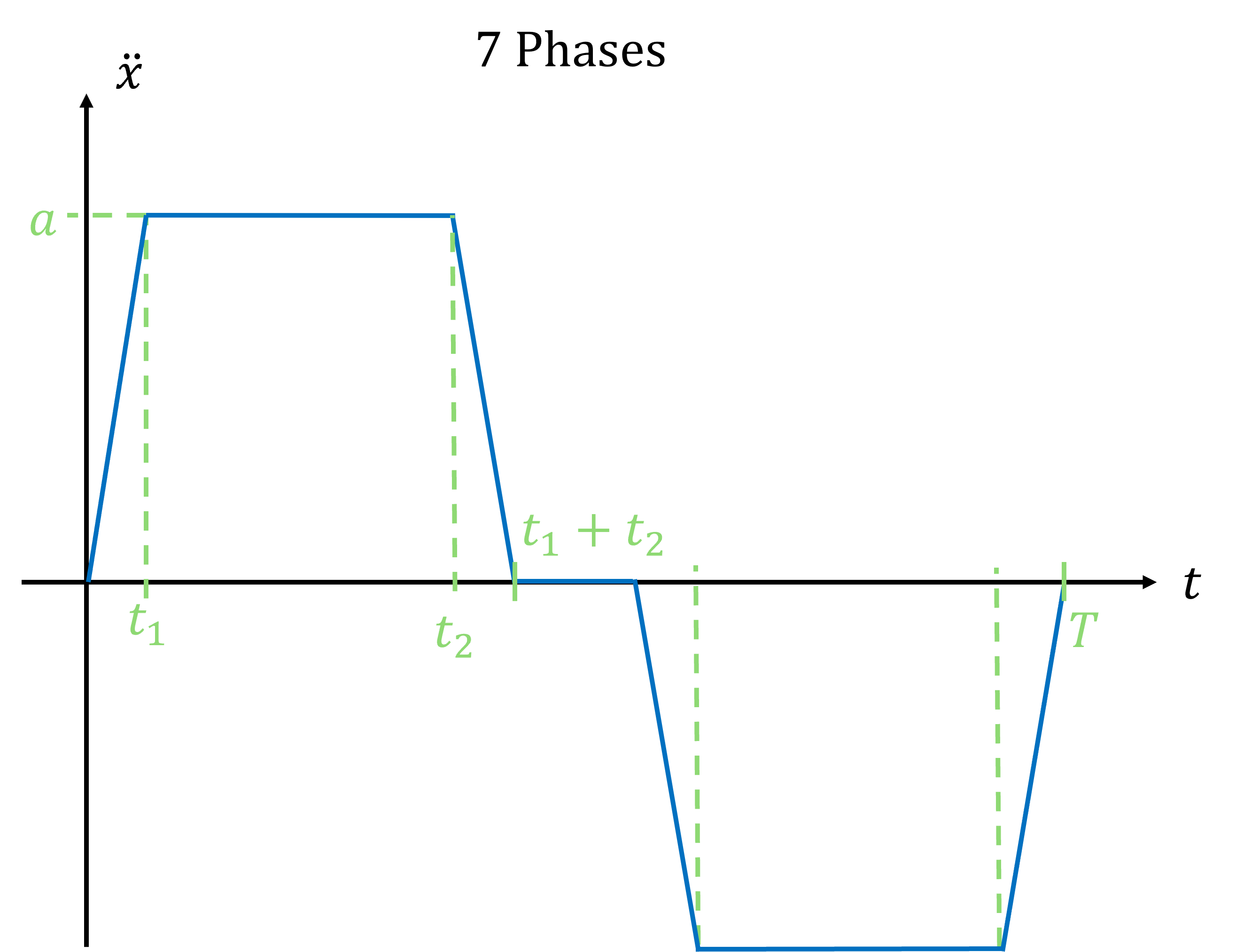}
	\caption{Acceleration profiles $\ddot{x}(t)$ for the four possible cases of a jerk-controlled motion.} \label{Figure2}
\end{figure}

By including the jerk constraint $(N=3)$ we now reach the actual case of the SPM. However, there are four cases to be distinguished (see Figure~\ref{Figure2} and Section~\ref{5.2}). We derive Equation~(\ref{Eins}) for each of the cases. Let the pre-optimal constant $c$ be denoted by $j$, let $v$ be defined as before in Section~\ref{3.1}, and let $a:=\max\limits_{t\in[0,T]}|\ddot{x}(t)|$ denote the maximal acceleration reached during the motion.\\

\noindent\textbf{4-Phase Case.} Clearly,
\begin{equation}
	a = jt_1 \quad \mathrm{and} \quad v = at_1 		\label{a=}
\end{equation}
and for $t\in[t_1,2t_1]$ we find 
	\begin{align} 
		\dddot{x}(t) &= -j \label{dddot}\\
		\ddot{x}(t) &= a - j  (t-t_1)\\
		\dot{x}(t) &= \frac{1}{2}v + a  (t-t_1) - \frac{1}{2}j  (t-t_1)^2\\
		x(t) &= \frac{1}{6}jt_1^3 + \frac{1}{2}v  (t-t_1) + \frac{1}{2}a  (t-t_1)^2 - \frac{1}{6}j(t-t_1)^3 \label{4-phase-case}
	\end{align} 
From this we get by $s = 2  x(2t_1)$ that $s = vt_1 + at_1^2$ or equivalently $\frac{s}{v} = t_1 + \frac{a}{v} t_1^2 = 2t_1$ with Equation~(\ref{a=}).
Finally,
\begin{equation}
	T = 4t_1 = 2t_1 + t_1+t_1 = \frac{s}{v} + \frac{v}{a} + \frac{a}{j}. 
\end{equation}

\noindent\textbf{5-Phase Case.} Again Equation~(\ref{a=}) is valid, but now we have with Equation~(\ref{4-phase-case}) that
\begin{equation}
	s = 2  x(2t_1) + v  (T-4t_1) = vt_1 + at_1^2 + v  (T-4t_1)
\end{equation}
which is equivalent to $\frac{s}{v} = T-2t_1$. We conclude via Equation~(\ref{a=}) that
\begin{equation}
	T = \frac{s}{v} + 2t_1 = \frac{s}{v} +t_1 + t_1 = \frac{s}{v} + \frac{v}{a} + \frac{a}{j}. 
\end{equation}
\noindent\textbf{6-Phase Case.} In this case, 
\begin{equation}
	a = jt_1 \quad \mathrm{but} \quad v = at_2 = a  \left( \frac{T}{2} - t_1\right).  \label{6-phase-case}
\end{equation}
Then for $t\in [t_1,t_2]$:
\begin{align}
	\ddot{x}(t) &= a\\
	\dot{x}(t) &= \frac{1}{2} at_1 + a  (t-t_1) \\
	x(t) &= \frac{1}{6} jt_1^3 + \frac{1}{2} at_1  (t-t_1) + \frac{1}{2} a  (t-t_1)^2
\end{align}
and thus $x(t_2) = \frac{1}{6} jt_1^3 + \frac{1}{2} at_1 (t_2-t_1) + \frac{1}{2} a (t_2-t_1)^2$. For $t \in [t_1,t_2]$ we have
\begin{align}
	\dddot{x}(t) &= -j \\
	\ddot{x}(t) &= a - j(t-t_2) \\
	\dot{x}(t) &= \frac{1}{2} at_1 + a (t_2-t_1) + a (t-t_2) - \frac{1}{2} j (t-t_2)^2 \\
	x(t) &= x(t_2) + \left(\frac{1}{2} at_1 + a (t_2-t_1) \right) (t_2-t_1) + \frac{1}{2} a (t-t_2)^2 - \frac{1}{6} j (t-t_2)^3
\end{align}
which leads -- after some manipulations -- to
\begin{equation}
	x(t_1+t_2) = x\left(\frac{T}{2} \right) = \frac{1}{2} a t_1t_2 +  \frac{1}{2} at_2^2 \label{xt_1}
\end{equation}
\begin{equation}
	s = 2  x\left(\frac{T}{2} \right) = at_1t_2 + at_2^2
\end{equation}
and via Equation~(\ref{6-phase-case}) to $\frac{s}{v} = \frac{T}{2}$. Applying Equation~(\ref{xt_1}) leads to $\frac{v}{a} + \frac{a}{j} = \frac{T}{2}$ as well, and hence
\begin{equation}
	T = \frac{T}{2} + \frac{T}{2} = \frac{s}{v} + \frac{v}{a} + \frac{a}{j}
\end{equation}
which again proves Equation~(\ref{Eins}).\\

\noindent\textbf{7-Phase Case.} From the 6-phase case we can already adapt
\begin{equation}
	a = jt_1, v=at_2 \quad \mathrm{and} \quad (t_1+t_2) \stackrel{(\ref{xt_1})}{=} \frac{1}{2} at_2(t_1+t_2) = \frac{1}{2} v(t_1+t_2) \label{7-phase-case}
\end{equation}
For the distance $s$ we compute
\begin{align}
	s &= 2 x (t_1+t_2) + v  (T-2(t_1+t_2)) \\
	&\stackrel{(\ref{xt_1})}{=} v  (t_1+t_2) + v  (T-2(t_1+t_2))
\end{align}
Thus, 
\begin{equation}
	\frac{s}{v} = T - t_1 - t_2 \stackrel{(\ref{7-phase-case})}{=} T - \frac{a}{j} - \frac{v}{a}
\end{equation}
which completes the proof of Equation~(\ref{Eins}) for all cases of SPM.

	\section{Time Horizon of Pre-Optimal Motions} \label{Section4}
However, since scanning larger decision trees seems quite unsatisfactory, this section presents a systematic approach for arbitrary orders.
\begin{theorem}\label{theo3}
Let $x(t)$ be a pre-optimal trajectory of $N$-th order. For $0 \leq n \leq N$ define
\begin{equation}
	x_n = \max\limits_{t\in[0,T]}|x^{(n)}(t)|.
\end{equation}
Then the time horizon $T$ of the trajectory is given by
\begin{equation}
	T = \sum_{n=0}^{N-1} \frac{x_n}{x_{n+1}}.
\end{equation}
\end{theorem}
Indeed, this is the desired generalization of Equation~(\ref{Eins}). However, it is possible to state an even stronger result which we will prove recursively. To this end, we need another definition:
\begin{definition}\label{Def4}
	With the notation from Theorem~\ref{theo3}, we introduce the \textup{earliest time point of achievement} as 
	\begin{equation}
		T_n = \min \left\lbrace t \in [0,T]~|~x^{(n)} (t) = x_n \right\rbrace  \label{Eq1Def4}
	\end{equation}
	for $0 \leq n \leq N-1$.
\end{definition}
Then the time horizon $T$ is equivalent to $T_0$. $T_1$ is exactly that time point where the velocity $x^{(1)} (t) = \dot{x}(t)$ initially reaches its maximum $x_1 = v$.
\begin{theorem}\label{theo5}
	Let $x(t)$ be a pre-optimal trajectory of order $N$. Then
	\begin{equation}
		T_n = \sum_{k=n}^{N-1}\frac{x_k}{x_{k+1}} \label{EqTheo5}
	\end{equation}
	for $0 \leq n \leq N-1$.
\end{theorem}
Note that Theorem~\ref{theo3} follows directly from Theorem~\ref{theo5} by $T=T_0$. To prove Theorem~\ref{theo5}, we will construct pre-optimal trajectories recursively starting at order $N$.
\begin{definition} \label{Def6}
	Consider a sequence $(t_n)_{n\geq0}$ with $t_0 > 0$ and $t_{n+1} \geq 2t_n$ for all $n$. Fix a constant $c_0 \in \mathbb{R}$ with $c_0 > 0$ and a function $f_0: [0,t_0] \rightarrow \mathbb{R}$ with $f_0(t) = c_0$ for all $t \in [0,t_0]$. Define the sequences of functions $(\tilde{f}_n)_{n\geq 0}, (f_n)_{n\geq 0}$ with $\tilde{f}_n, f_n: [0,t_n] \rightarrow\mathbb{R}$ and
	\begin{eqnarray}
	\tilde{f}_{n+1}(t)&=&	
	\begin{cases}
			f_n(t)	&\mbox{for} \quad 0\leq t\leq t_n\\
			-f_n(t_{n+1}-t)	&\mbox{for} \quad t_{n+1}-t_n \leq t \leq t_{n+1}\\
			0	&\mbox{elsewhere}
		\end{cases} \\
		f_{n+1}(t) &=& \int\limits_{0}^{t} \tilde{f}_{n+1}(\tilde{t}) \mathrm{d}\tilde{t}
	\end{eqnarray}
\end{definition}
The next lemma collects some basic properties of the sequence $(f_n)_{n\geq 0}$.

\begin{lemma}\label{Lem7}
	One has
\begin{enumerate}[label=(\roman*)]
	\item $f_n(t) \geq 0$ for all $t\in[0,t_n]$ and $n \geq 0$. \label{L7,i}
	\item $f_n(t) = f_n(t_n-t)$ for all $t\in[0,t_n]$ and $n \geq 0$. \label{L7,ii}
	\item $f_n(0) = f_n (t_n) = 0$ for all $n \geq 1$. \label{L7,iii}
	\item $f_n(t) = 0$ for all $t\in \left]t_{n-1}, t_n - t_{n-1} \right[$ and $n\geq 1$. \label{L7,iv}
	\item $\max\limits_{t\in[0,t_n]} f_n(t) = f_n(t_{n-1})$ for all $n\geq1$. \label{L7,v}
	\item $f_n(t) = f_n(t_{n-1}) - f_n(t_{n-1} - t)$ for all $0 \leq t \leq t_{n-1}$ and $n\geq1$. \label{L7,vi}
\end{enumerate}
\end{lemma}
\begin{proof}
	The proof of all statements is mainly based on the fact that $\dot{f}_n(t) = \tilde{f}_n(t) = f_{n-1}(t)$ for $n>0$, and the symmetry property $\tilde{f}_n(t) = -\tilde{f}_n (t_n - t)$. In more detail:
	\begin{enumerate}[label=(\roman*)]
		\item By induction: The case $n=0$ is clear due to the definition. Consider now $f_{n+1}$ while assuming $f_n(t) \geq 0$ for $t\in[0,t_n]$. Then $f_{n+1} (t) \geq 0$ for $t\in[0,t_n]$ since it is the integral of a positive function. For $t\in \left] t_n, t_{n+1} - t_n \right[, f_{n+1}(t)$ remains constant (and positive) because the integrand is $0$. For $t\in[t_{n+1} - t_n, t_{n+1}]$, the absolute value of the negative integral is at most as large as its positive counterpart. Hence, $f_{n+1}(t) \geq 0$ for all $t\in[0, t_{n+1}]$.
		\item This is an immediate consequence of the symmetry property.\label{Proofii}
		\item This follows from~\ref{Proofii} for $t=0$ and the definition.
		\item This follows from the fact that $\tilde{f}_n(t) = 0$ for $t_{n-1} < t < t_1 - t_{n-1}$.
		\item This is a consequence of the fact that $f_n$ is monotonously increasing for $0 \leq t \leq t_{n-1}$ since it has a positive derivative: $\dot{f}_n(t) = \tilde{f}_n(t) = f_{n-1}(t) \geq 0$ per definition.
		\item Integration of~\ref{Proofii} yields:
		\begin{align}
			\int\limits_{0}^{t} f_{n-1}(\tilde{t}) \mathrm{d}\tilde{t} &= \int\limits_{0}^{t} f_{n-1} (t_{n-1} - \tilde{t})\mathrm{d}\tilde{t}\\
			f_n(t) - f_n(0) &= (-1)(f_{n-1}(t_{n-1} - t) - f_n (t_{n-1}))
		\end{align}
		which gives the assertion after inserting $f_n(0) = 0$.
	\end{enumerate}
\end{proof}
	A direct consequence of Lemma~\ref{Lem7} is the following Lemma~\ref{Lem8}.
	\begin{lemma}\label{Lem8}
		For all $n \geq 1$ one has
		\begin{equation}
			\frac{f_{n+1}(t_n)}{f_n(t_{n-1})} = t_n - t_{n-1}.
		\end{equation}
	\end{lemma}
\begin{proof}
	Consider
	\begin{eqnarray}
	&&	\int\limits_{0}^{t_{n-1}} f_n(t) - \frac{f_n(t_{n-1})}{t_{n-1}} t \mathrm{d}t \nonumber\\ &=& \int\limits_0^{\frac{t_{n-1}}{2}} f_n(t) - \frac{f_n(t_{n-1})}{t_{n-1}} t \mathrm{d}t + \int\limits_{\frac{t_{n-1}}{2}}^{t_{n-1}} f_n(t) - \frac{f_n(t_{n-1})}{t_{n-1}} t \mathrm{d}t \nonumber\\
		&\overset{\mathrm{Lemma~}\ref{Lem7},~\ref{L7,vi}}{=}& \int\limits_{0}^{t_{\frac{n-1}{2}}} f_n(t) - \frac{f_n(t_{n-1})}{t_{n-1}} t \mathrm{d}t + \int\limits_{\frac{t_{n-1}}{2}}^{t_{n-1}} f_n(t_{n-1}) - f_n (t_{n-1}-t) - \frac{f_n(t_{n-1})}{t_{n-1}} t \mathrm{d}t \nonumber\\
		&\overset{\tilde{t}:=t_{n-1} - t}{=}& \int\limits_0^{\frac{t_{n-1}}{2}} f_n(t) - \frac{f_n(t_{n-1})}{t_{n-1}} t \mathrm{d}t + \int\limits_{\frac{t_{n-1}}{2}}^0 f_n(t_{n-1}) - f_n(\tilde{t}) - \frac{f_n(t_{n-1})}{t_{n-1}} \mathrm{d}\tilde{t} \nonumber\\
		&=& \int\limits_0^{\frac{t_{n-1}}{2}} f_n(t) - \frac{f_n(t_{n-1})}{t_{n-1}} t \mathrm{d}t - \int\limits_{0}^{\frac{t_{n-1}}{2}} f_n(\tilde{t}) - \frac{f_n(t_{n-1})}{t_{n-1}} \tilde{t} \mathrm{d}\tilde{t}\nonumber\\
		&=& 0
	\end{eqnarray}
	from which we conclude that
	\begin{equation}
		\int\limits_{0}^{t_{n-1}} f_n(t) \mathrm{d}t
 = \int\limits_{0}^{t_{n-1}} \frac{f_n(t_{n-1})}{t_{n-1}} t\mathrm{d}t = \frac{1}{2} t_{n-1} f_n(t_{n-1}).
 	\end{equation}
 	Due to Lemma~\ref{Lem7}~\ref{L7,ii} and \ref{L7,iv} we compute on the one hand 
 \begin{eqnarray}
 	\int\limits_{0}^{t_n} f_n(t) \mathrm{d}t &=& 2 \int\limits_{0}^{t_{n-1}}f_n(t) \mathrm{d}t + (t_n - 2t_{n-1}) f_n(t_{n-1})\nonumber\\
 	&=& t_{n-1} f_n(t_{n-1}) + (t_n - 2t_{n-1}) f_n(t_{n-1})\nonumber\\
 	&=& (t_n - t_{n-1}) f_n(t_{n-1})
 \end{eqnarray}
 and on the other hand
 \begin{equation}
 	f_{n+1}(t_n) \overset{\mathrm{Def.}~\ref{Def6}}{=} \int\limits_0^{t_n} \tilde{f}_{n+1}(t) \mathrm{d}t \overset{\mathrm{Def.}~\ref{Def6}}{=}\int\limits_0^{t_n} f_n(t) \mathrm{d}t = (t_n - t_{n-1}) f_n(t_{n-1}) 
 \end{equation}
 which proves the assertion.
\end{proof}
\begin{proof}[Proof of Theorem~\ref{theo5}.]
	For every pre-optimal trajectory $x(t)$ of order $N$ there exist sequences $(t_n)_{n\geq0}$, $(f_n)_{n\geq0}$ according to Definition~\ref{Def6} such that $x(t) = f_N(t)$ for $0 \leq t \leq t_{N-1} = T$. Then $x^{(n)}(t) = f_{N-n}(t)$ for $0 \leq n \leq N$, and in addition
	\begin{eqnarray}
		T_n &=& \min\left\lbrace t \in[0,T]~|~x^{(n)}(t) = x_n\right\rbrace \nonumber\\
		&=& \min\left\lbrace t \in[0,T]~|~f_{N-n}(t) = f_{N-n}(t_{N-n-1})\right\rbrace \nonumber\\
		&=& t_{N-n-1}
	\end{eqnarray}
	With Lemma~\ref{Lem8} we arrive at
	\begin{equation}
		T_n - T_{n+1} = \frac{x_n}{x_{n+1}} \label{Arrive}
	\end{equation}
	and regarding $n=N$ we have $T_N = \frac{x_{N-1}}{x_N}$ as starting point of the recursion~(\ref{Arrive}). Finally, a straight forward calculation shows that $T_n$ given by Equation~(\ref{EqTheo5}) is the explicit solution of Equation~(\ref{Arrive}).
\end{proof}
	\section{Applications and Examples} \label{Section5}
Equipped with the general Equation~(\ref{EqTheo5}) of Theorem~\ref{theo5} we are able to gain some new insights into of optimal motion planning. In the first subsection we will briefly discuss the cases $N=1$ (velocity-controlled motions) and $N=2$ (acceleration-controlled motions) before the might and elegance of Equation~(\ref{Eins}) are fully exploited while treating jerk-controlled and snap-controlled movements (Subsections~\ref{5.2} and~\ref{5.3}, respectively). Section~\ref{Section5} closes with remarks on arbitrary order (Subsection~\ref{5.4}) which culminate in the algorithm for calculating the minimum-time motion (i.e. MTT, Subsection~\ref{5.5}).
\subsection{Velocity-Controlled and Acceleration-Controlled Motions}\label{5.1}
As mentioned earlier (see beginning of Section~\ref{Section3}), the case $N=1$, which is a uniform motion, is rather simple since the only relevant equation is
\begin{equation}
	T = \frac{s}{v}. \label{Eq5.1}
\end{equation}
Recall that $v=x_1=\max\limits_{t\in[0,T]}|\dot{x}(t)|$ denotes the maximal reached velocity during the movement. Regarding MTT, the obvious optimum is achieved if $v=w_1=v_{\max}$ (velocity bound, see Definition~\ref{Def1}). Then $T_{\mathrm{MTT}} = \frac{s}{v_{\max}}$. Conversely, if both distance $s$ and a time horizon $T_{\mathrm{MDT}}$ are given (clearly, $T_{\mathrm{MDT}} \geq T_{\mathrm{MTT}}$), then the solution of MDT (for $M=1$) is the uniform motion with the minimal velocity given by $v_{\mathrm{MDT}}= \frac{s}{T_{\mathrm{MDT}}}\leq v_{\max}$. Nevertheless, Equation~(\ref{Eq5.1}) is the very inconspicuous sister of Equation~(\ref{Eins}).

\begin{figure}[h!]
	\centering
	\includegraphics[width=0.5\textwidth]{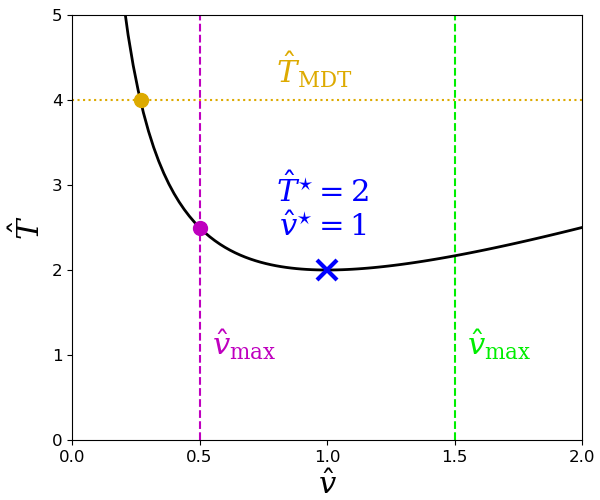}
	\caption{$\hat{T}$ as a function of $\hat{v}$ (black curve) according to Equation~(\ref{hatT}) which exhibits a global minimum at $\hat{v}^{\star} = 1$ and  $\hat{T}^{\star} = 2$ (blue cross). The triangular case of MTT occurs if $\hat{v}_{\max}$ lies beyond this minimum (green dashed line). Then the minimal time is exactly $\hat{T}^{\star}$. Otherwise (magenta dashed line, trapezoidal case) the optimum of MTT is located at the magenta bullet. The MDT optimization is colored in orange. The optimum is indicated by the orange bullet.} \label{Figure3}
\end{figure}

The situation becomes a bit more interesting for $N=2$ (acceleration control). Let $v$ be defined as above, set $a=x_2=\max\limits_{t\in[0,T]}|\ddot{x}(t)|$, and $v_{\max} =w_1$ as well as $a_{\max} = w_2$. Due to Theorem~\ref{theo5} we have
\begin{align}
	T_0 &= T = \frac{s}{v} + \frac{v}{a} \label{T_0}\\
	T_1 &= \frac{v}{a}
\end{align}
where $T_1$ is the time after which $v$ is initially achieved. Obviously, the solutions of MTT and MDT ($M=1$) both obey $a=a_{\max}$. Hence we can consider a reduced problem by introducing dimensionless quantities (denoted with $\hat{}$~):
\begin{equation}
	\hat{T} = T a_{\max}^{\frac{1}{2}} s^{-\frac{1}{2}} \qquad
	\hat{v} =  va_{\max}^{-\frac{1}{2}}  s^{-\frac{1}{2}} \qquad
	\hat{s} = \hat{a}_{\max} = 1
\end{equation}
which leads to
\begin{equation}
	\hat{T} = \frac{1}{\hat{v}} + \hat{v} \label{hatT}
\end{equation}
and $\hat{v}_{\max}$ is transformed analogously. Figure~\ref{Figure3} reveals the situation: $\hat{T}(\hat{v})$ has a global minimum at $\hat{v}^{\star} = 1$ and $\hat{T}^{\star} = 2$, and the bound $\hat{v}_{\max}$ is a vertical cut-off. Thus, the solution of MTT is either
\begin{itemize}
	\item $\hat{v}_{\mathrm{MTT}} = \hat{v}^{\star}$, $\hat{T}_{\mathrm{MTT}} = \hat{T}^{\star}$, if $\hat{v}^{\star} \leq \hat{v}_{\max}$ (triangular case in Section~\ref{3.1}, see Figure~\ref{Figure1} left), or 
	\item  $\hat{v}_{\mathrm{MTT}} = \hat{v}_{\max}$, $\hat{T}_{\mathrm{MTT}} = \frac{1}{\hat{v}_{\max}} + \hat{v}_{\max}$, otherwise (trapezoidal case, see Figure~\ref{Figure1} right).
\end{itemize}
Regarding MDT with a given $\hat{T}_{\mathrm{MDT}} \geq \hat{T}_{\mathrm{MTT}}$ (a horizontal cut-off), the minimal velocity $\hat{v}_{\mathrm{MDT}}$ is given by the smallest solution of the algebraic equation
\begin{equation}
	\hat{T}_{\mathrm{MDT}} = \frac{1}{\hat{v}_{\mathrm{MDT}}} + \hat{v}_{\mathrm{MDT}}.
\end{equation}
As a final remark, note that the global minimum can be also obtained by solving the equation
\begin{equation}
	\hat{T}_0 = 2\hat{T}_1 \label{final5.1}
\end{equation}
which will be revisited later.
\subsection{Jerk-Controlled Motions}\label{5.2}
According to Theorem~\ref{theo5} we have
\begin{align}
	T_0 &= T = \frac{s}{v} + \frac{v}{a} + \frac{a}{j} \label{5.2T0}\\
	T_1 &= \frac{v}{a} + \frac{a}{j} \label{5.2T1}\\
	T_2 &= \frac{a}{j} \label{5.2T2}
\end{align}
with the notation as above. Additionally, set $j_{\max} = w_3$ and $j = x_3 = \max\limits_{t\in[0,T]} |\dddot{x}(t)|$ (recall Definition~\ref{Def1}). Recall the systematics of the time variables: $T_0$ is the time after which the distance is accomplished; $T_1$ and $T_2$ are the time intervals after which $v$ and $a$ are initially achieved, respectively. Again, we perform a transformation to dimensionless quantities:
\begin{equation}
	\hat{T} = Tj_{\max}^{\frac{1}{3}}s^{-\frac{1}{3}} \qquad \hat{v} = vj_{\max}^{-\frac{1}{3}}s^{-\frac{1}{3}} \qquad \hat{a} = aj_{\max}^{-\frac{2}{3}}s^{-\frac{1}{3}} \qquad \hat{\jmath} = \hat{s} = 1 \label{5.2hatT}
\end{equation}
to arrive at
\begin{align}
	\hat{T}_0 &= \hat{T} = \frac{1}{\hat{v}} + \frac{\hat{v}}{\hat{a}} + \hat{a} \label{5.2hatT0}\\
	\hat{T}_1 &= \frac{\hat{v}}{\hat{a}} + \hat{a} \label{hatT1}\\
	\hat{T}_2 &= \hat{a}. \label{hatT2}
\end{align}
The parameter $j_{\max}$ is used in Equation~(\ref{5.2hatT}) because for all solutions of MTT and MDT ($M=1$ and $M=2$) it is known that $j = j_{\max}$ (cf. \cite{R1}).

\begin{figure}[h!]
	\centering
	\includegraphics[width=0.5\textwidth]{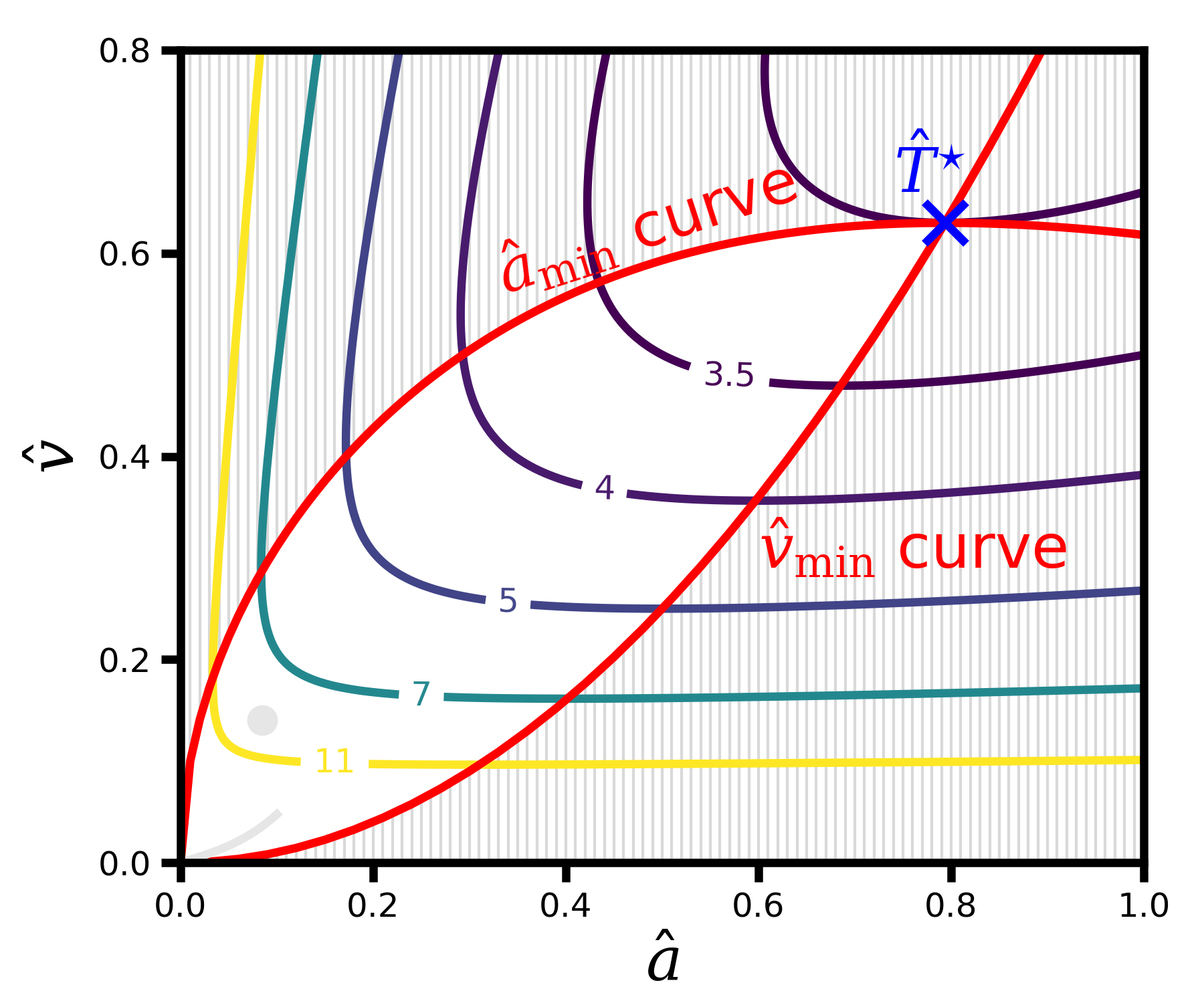}
	\caption{Contour plot of $\hat{T}$ as a function of $\hat{a}$ and $\hat{v}$ according to Equation~(\ref{5.2hatT0}). The global minimum (blue cross) is located at $\hat{a}_{\star} = \sqrt[3]{\frac{1}{2}}$, $\hat{v}^{\star} = \sqrt[3]{\frac{1}{4}}$ and $\hat{T}^{\star} = \sqrt[3]{32}$. The red solid boundary curves of the admissible region are given by $\hat{v} = \hat{a}^2$ ($\hat{v}_{\min}$ curve) and $\hat{v} = -\frac{1}{2} \hat{a}^2+ \sqrt{\frac{1}{4} \hat{a}^4 + \hat{a}}$ ($\hat{a}_{\min}$ curve).} \label{Figure4}
\end{figure}

The ``diving fish diagram'' of Figure~\ref{Figure4} displays $\hat{T}$ as a function of $\hat{v}$ and $\hat{a}$ as a contour plot. There are two red curves -- $\hat{a}_{\min}$ curve and $\hat{v}_{\min}$ curve -- which delimit the region of all principally possible motions. These curves are defined by the conditions (see Equation~(\ref{dddot})):
\begin{align}
	\hat{T}_0 &= 2\hat{T}_1 \quad (\hat{a}_{\min}~\mathrm{curve}) \label{amin}\\
	\hat{T}_1 &= 2\hat{T}_2 \quad (\hat{v}_{\min}~\mathrm{curve}) \label{vmin}
\end{align}
If such a condition is satisfied the associated order has no cruising phase. For example, $\hat{T}_1 = 2\hat{T}_2$ means that $\ddot{x}$ has a triangular shape if it is not equal to $0$ (see Figure~\ref{Figure2} both upper panels). If there is no cruising phase in $\ddot{x}$, the time for reaching $\hat{v}$ (equal to $\hat{T}_1$) is exactly twice the time needed for reaching $\hat{a}$ (equal to $\hat{T}_2$). This corresponds to the fact that for every acceleration a certain velocity is quasi automatically achieved. Therefore, not all of the ($\hat{v}, \hat{a}$)-pairs are actually possible. And it is precisely these borderline cases that are reflected in the conditions~(\ref{amin},~\ref{vmin}).

Furthermore, the global minimum is determined by simultaneously satisfying the equations~(\ref{amin},~\ref{vmin}). It is located at $\hat{a}^{\star} = \sqrt[3]{\frac{1}{2}}$, $\hat{v}^{\star} = \sqrt[3]{\frac{1}{4}}$ and $\hat{T}^{\star} = \sqrt[3]{32}$.

\begin{figure}[h!]
	\centering
	\includegraphics[width=0.5\textwidth]{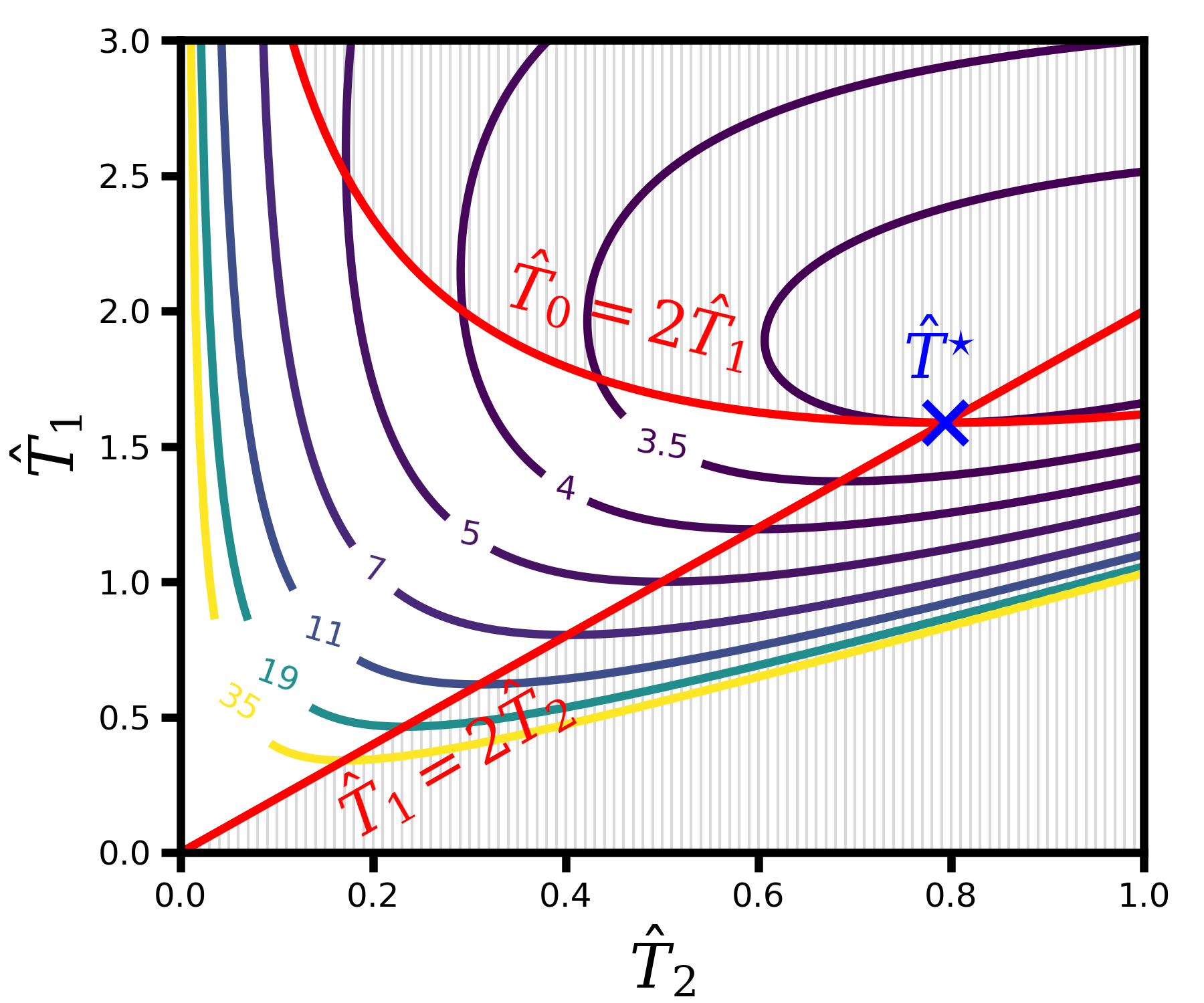}
	\caption{Dimensionless time horizon $\hat{T}$ as a function of $\hat{T}_1$ (dimensionless time for achieving $\hat{v}$) and $\hat{T}_2$ (respective time for achieving $\hat{a}$) which is basically a transformation of the diving fish diagram displayed in Figure~\ref{Figure4}. Note that $\hat{T} = ((\hat{T}_1 - \hat{T}_2)\hat{T}_2)^{-1} + \hat{T}_1$. The red straight line indicated by $\hat{T}_1 = 2\hat{T}_2$ corresponds to the $\hat{v}_{\min}$ curve (see Equation~(\ref{vmin})), and the red curve indicated by $\hat{T}_0=2\hat{T}_1$ belongs to the $\hat{a}_{\min}$ curve (see Equation~(\ref{amin})). This curve is described by $\hat{T}_1 = \frac{1}{2}\hat{T}_2+ \sqrt{\frac{1}{4}\hat{T}_2^2 + \hat{T}_2^{-1}}$.} \label{Figure4B}
\end{figure}

\begin{figure}[h!]
	\centering
	\includegraphics[width=0.4\textwidth]{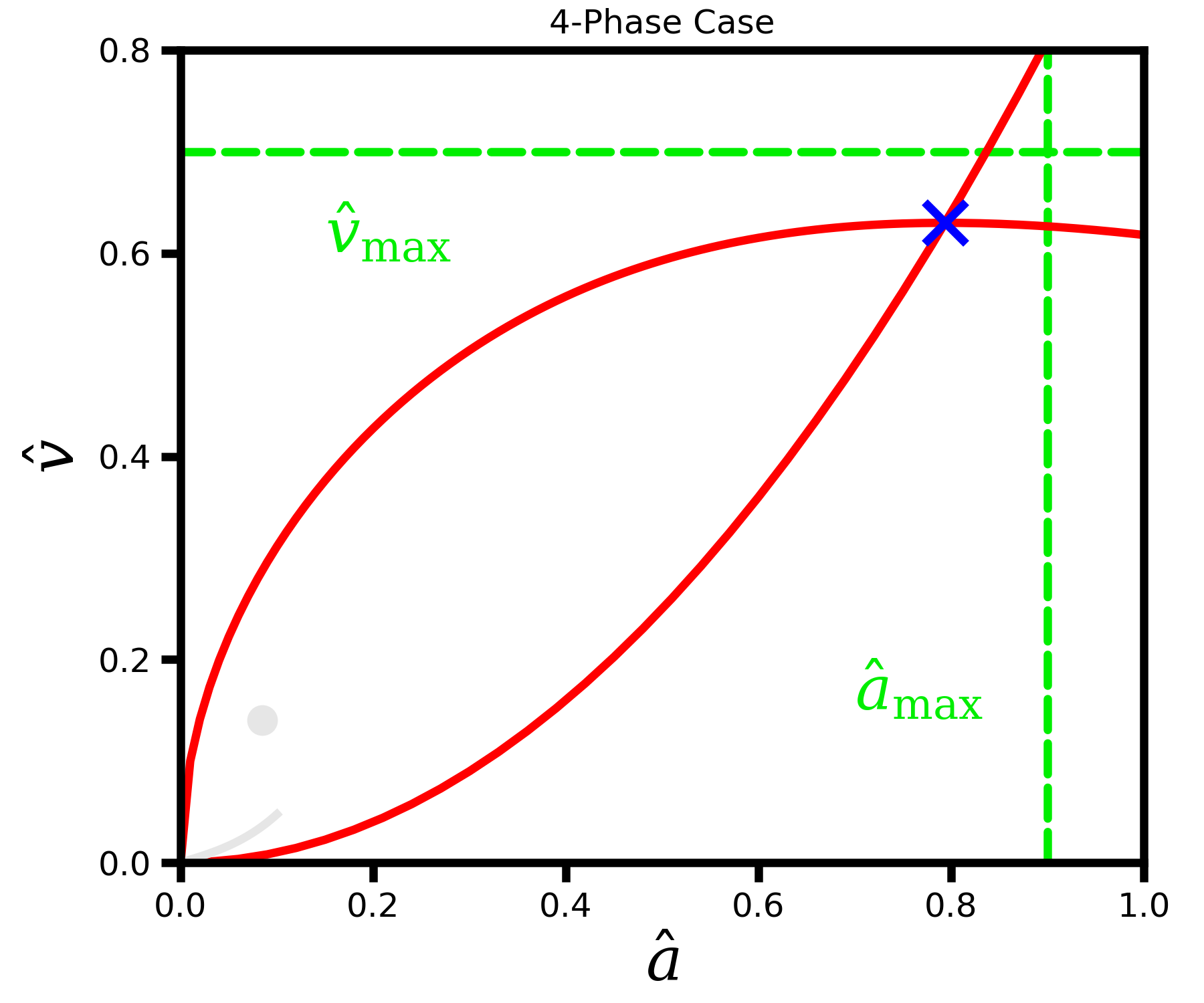}
	\includegraphics[width=0.4\textwidth]{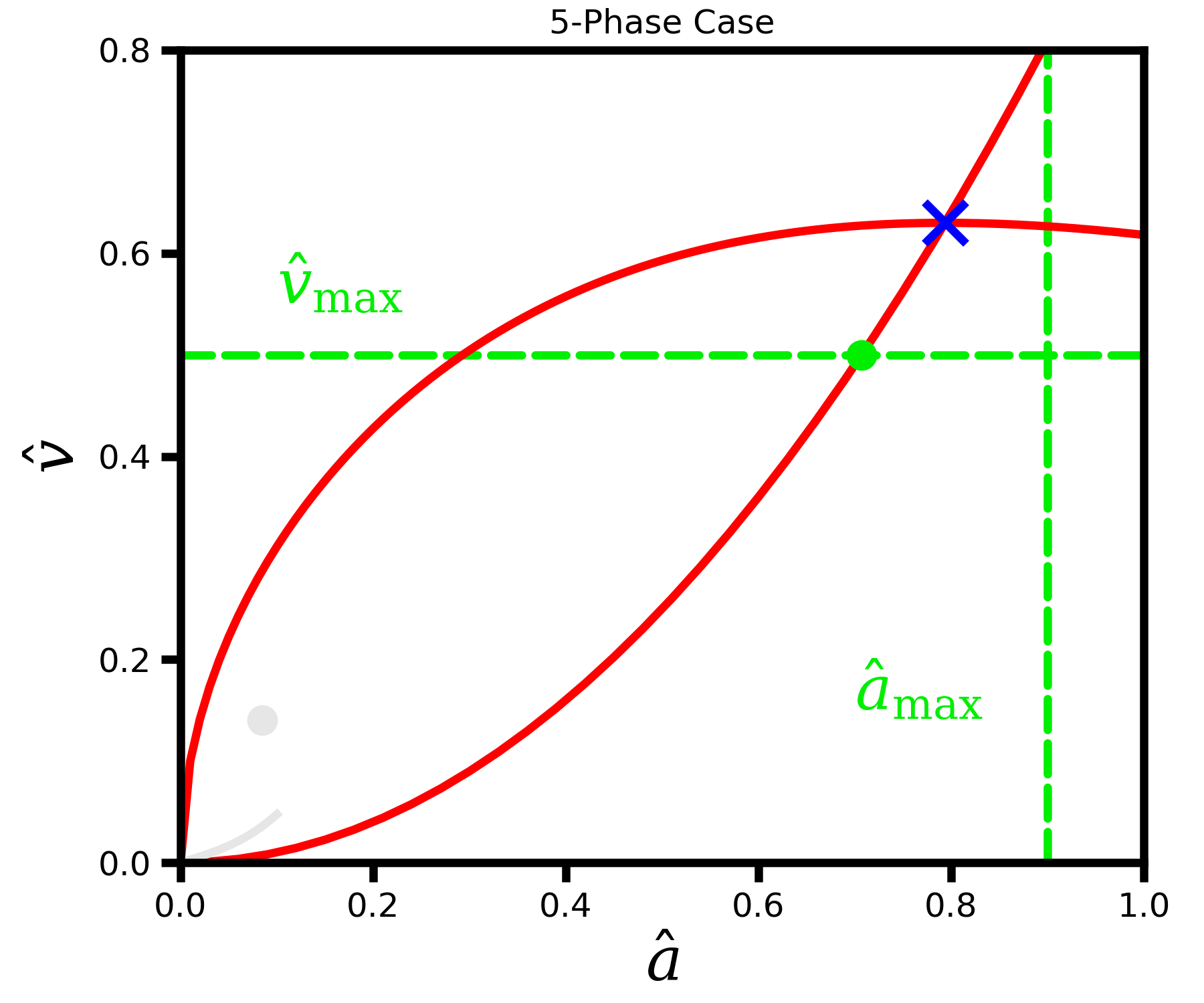}\\
	\includegraphics[width=0.4\textwidth]{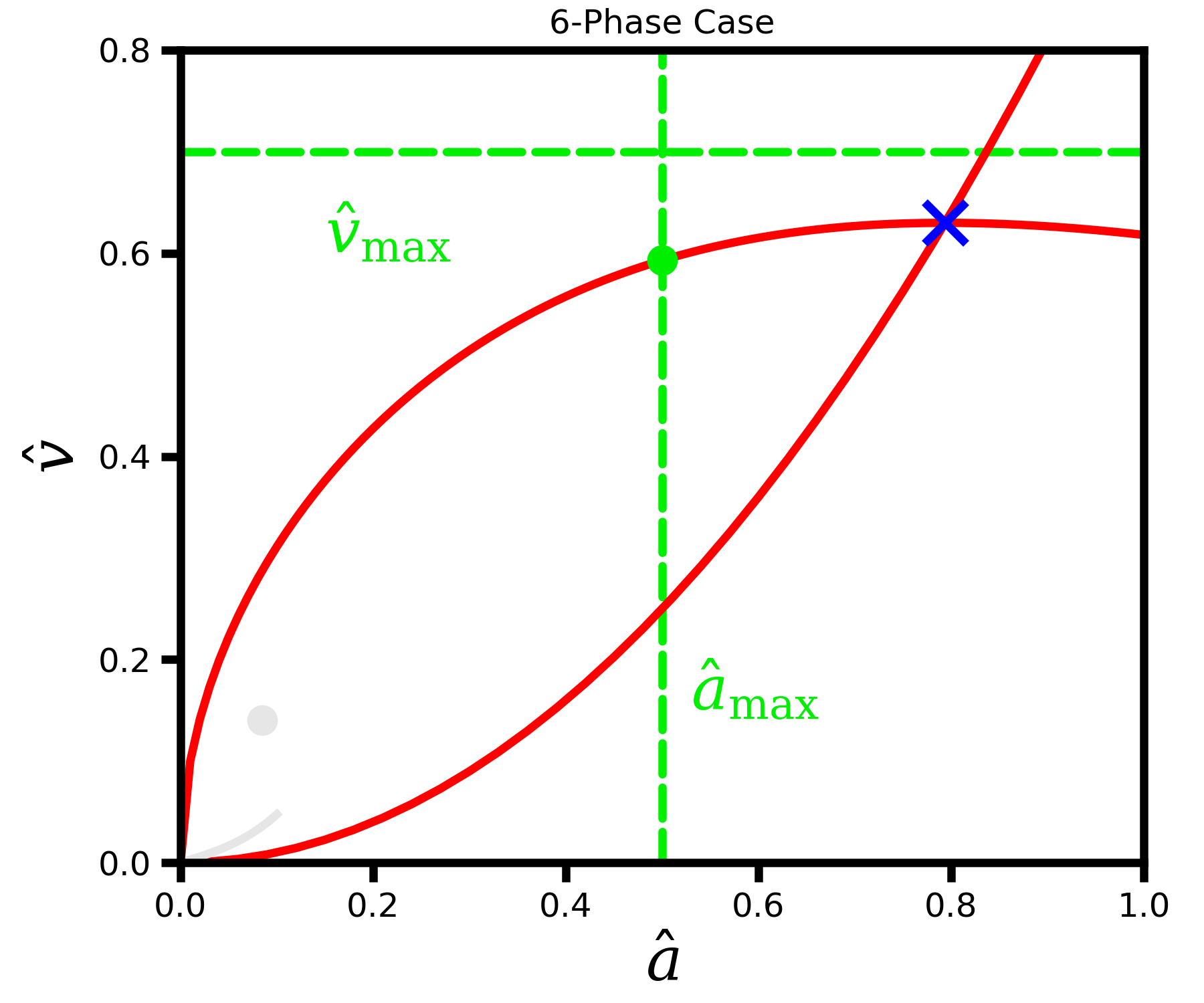}
	\includegraphics[width=0.4\textwidth]{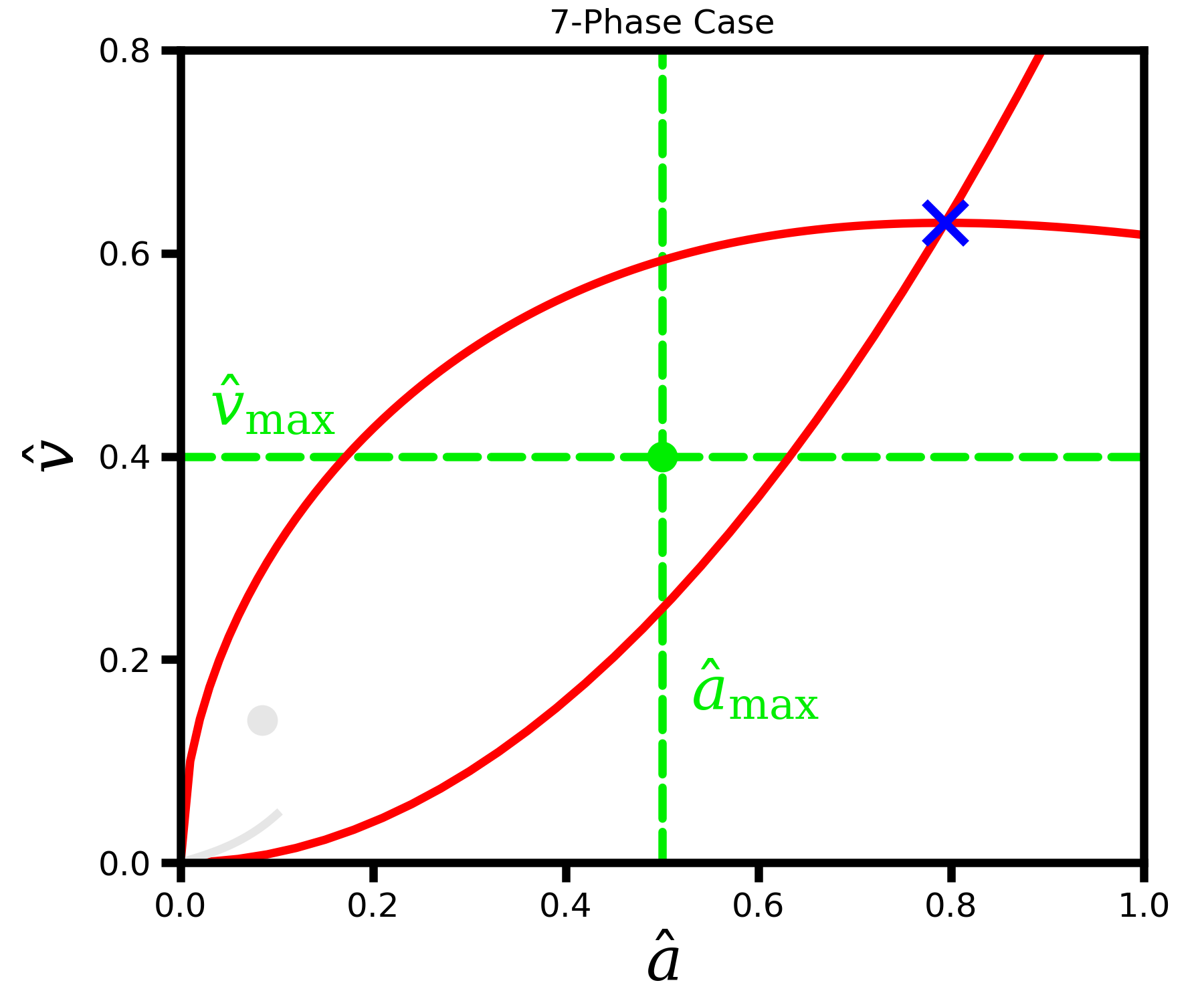}
	\caption{Schematic overview of the 4 possible cases of a jerk-controlled motion distinguished by the number of phases. The solution of MTT is marked with a green bullet (except for the 4-phase case, where the solution is equivalent to the global minimum highlighted by the blue cross). In each case, the admissible region lies below and left the $\hat{v}_{\max}$ and $\hat{a}_{\max}$ bounds, respectively, and is additionally bounded by the red curves.} \label{Figure5}
\end{figure}

\begin{figure}[h!]
	\centering
	\includegraphics[width=0.4\textwidth]{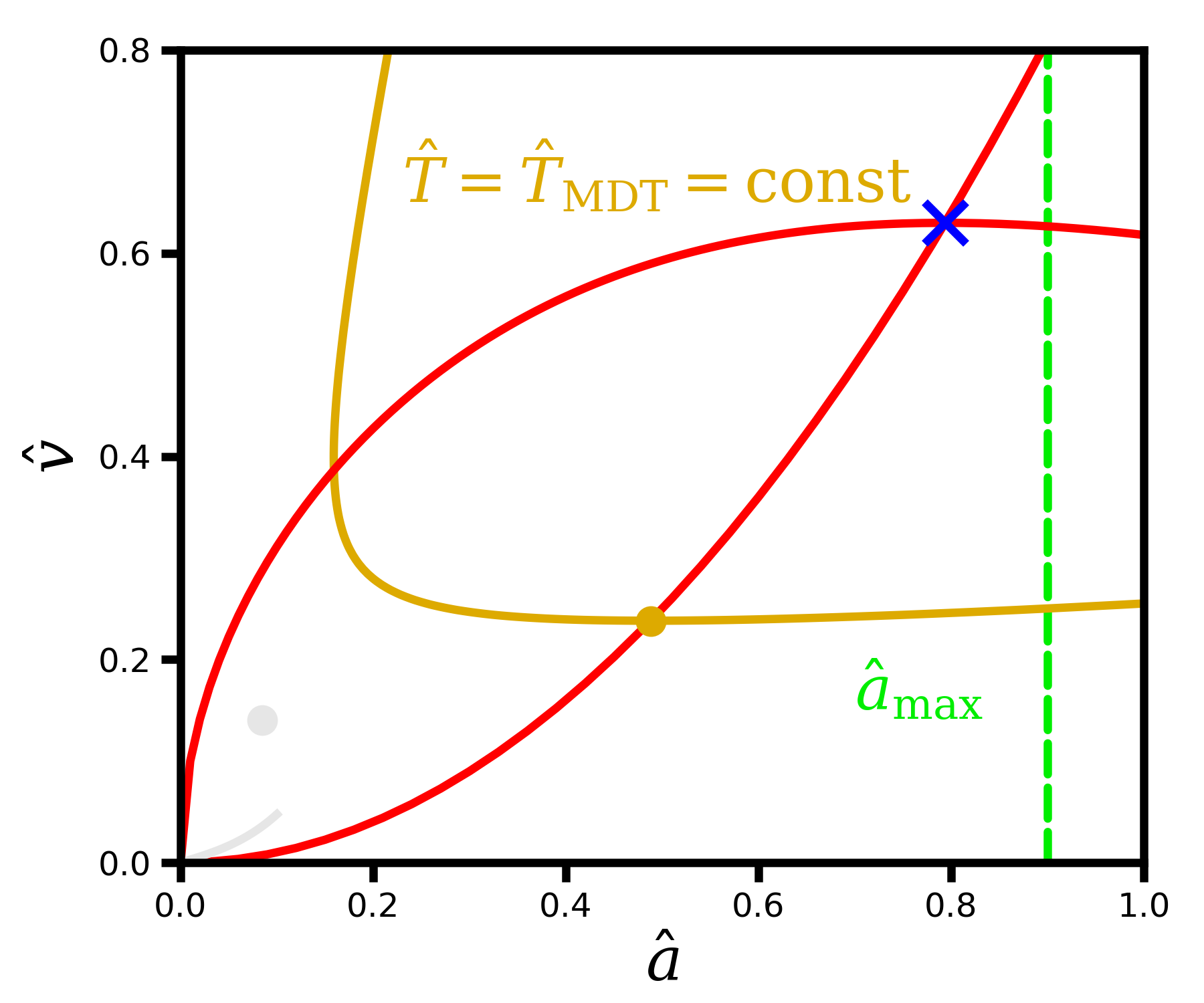}
		\includegraphics[width=0.4\textwidth]{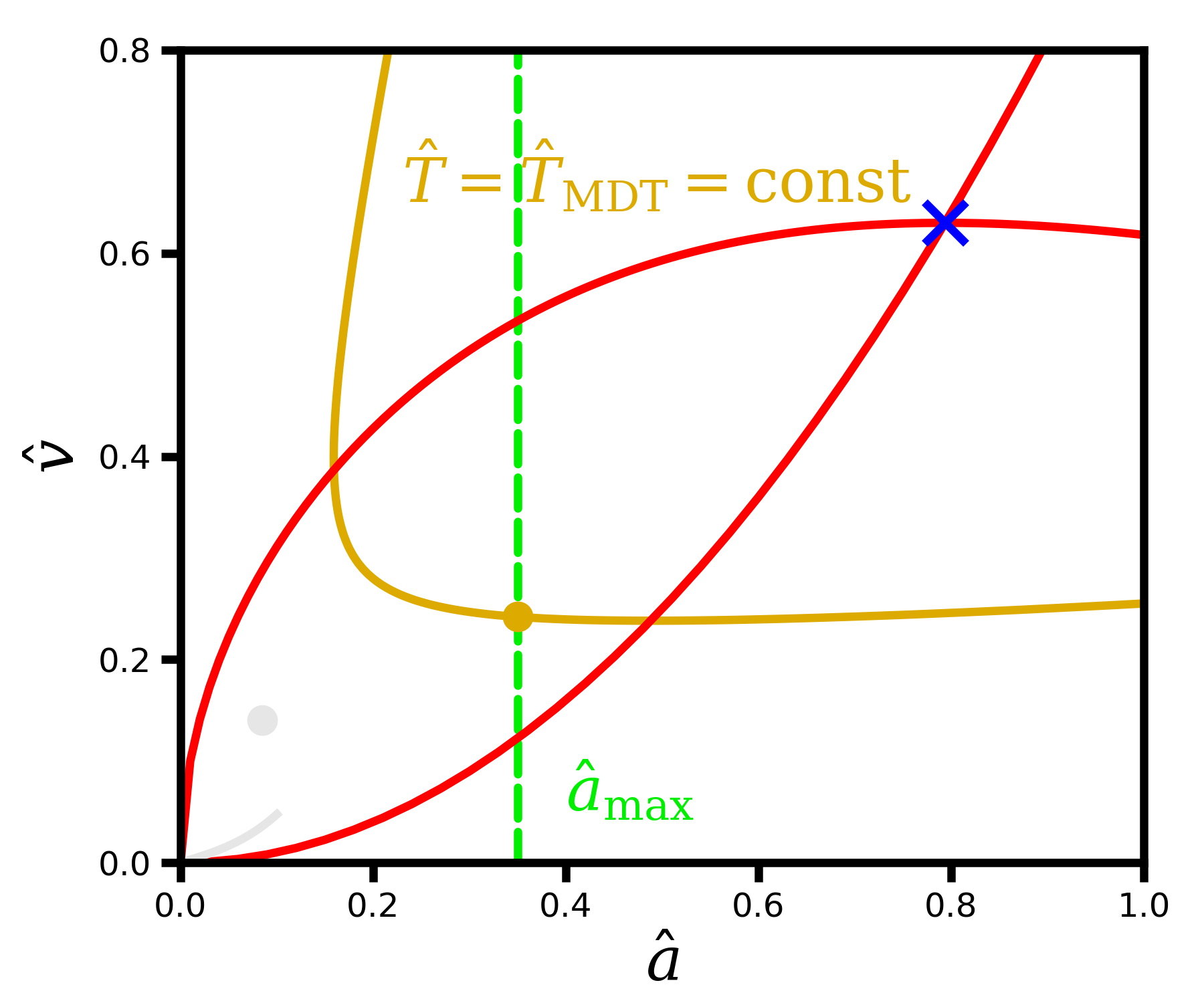}
	\caption{Schematic illustration of the solution of MDT with $M=1$, i.e. minimum-velocity motions. In any case, the solution is indicated by an orange bullet. Left: if $\hat{a}_{\max} > \hat{a}^{\star}$, the minimal velocity is achieved at the point where the contour curve belonging to the given time horizon $\hat{T}_{\mathrm{MDT}}$ leaves the admissible region and intersects the $\hat{v}_{\max}$ curve. Right: For small enough values of $\hat{a}_{\max}$ the contour curve reaches the acceleration boundary before the $\hat{v}_{\min}$ curve. The solution is determined by this intersection.} \label{Figure6}
\end{figure}

\begin{figure}[h!]
	\centering
	\includegraphics[width=0.4\textwidth]{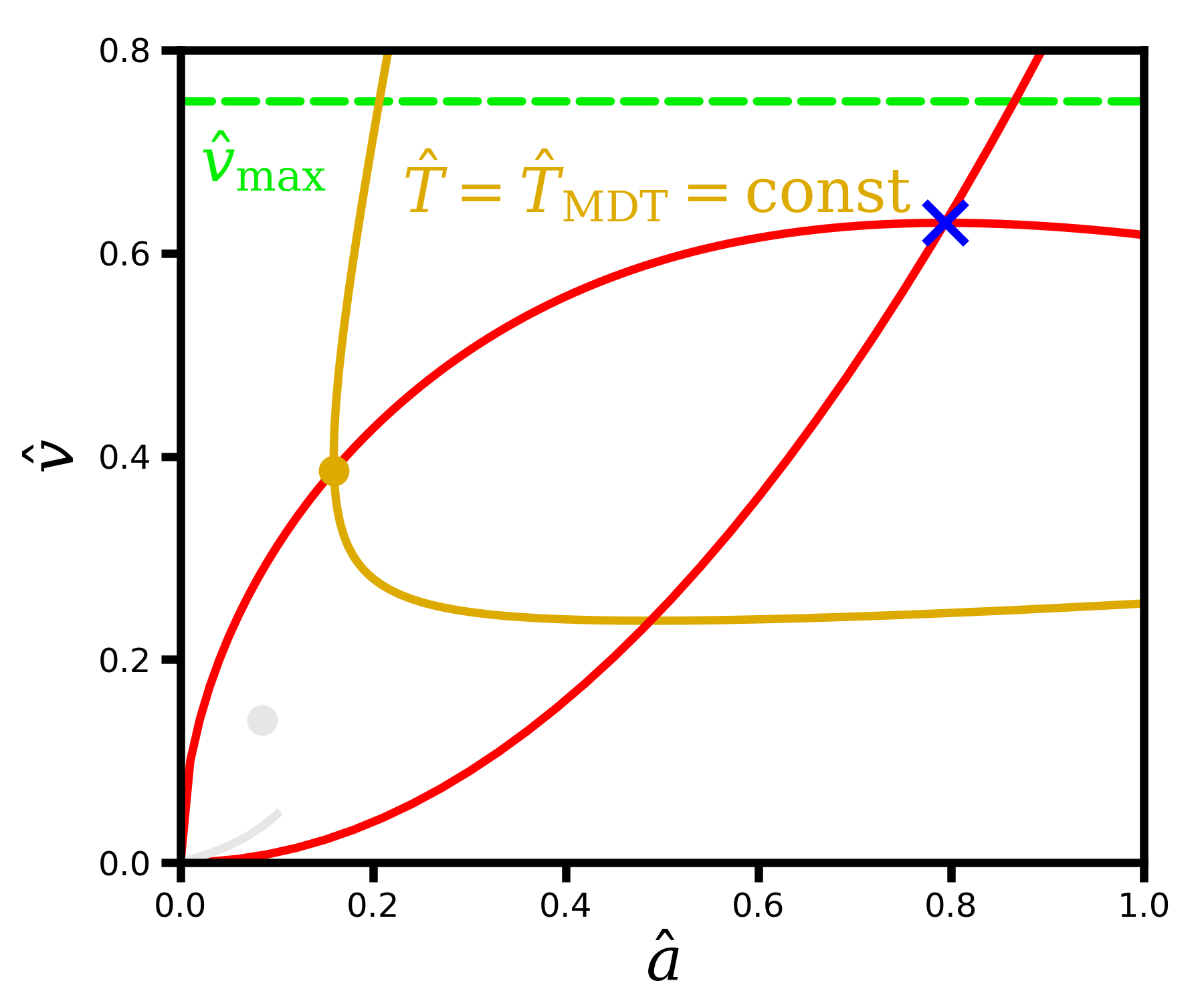}
	\includegraphics[width=0.4\textwidth]{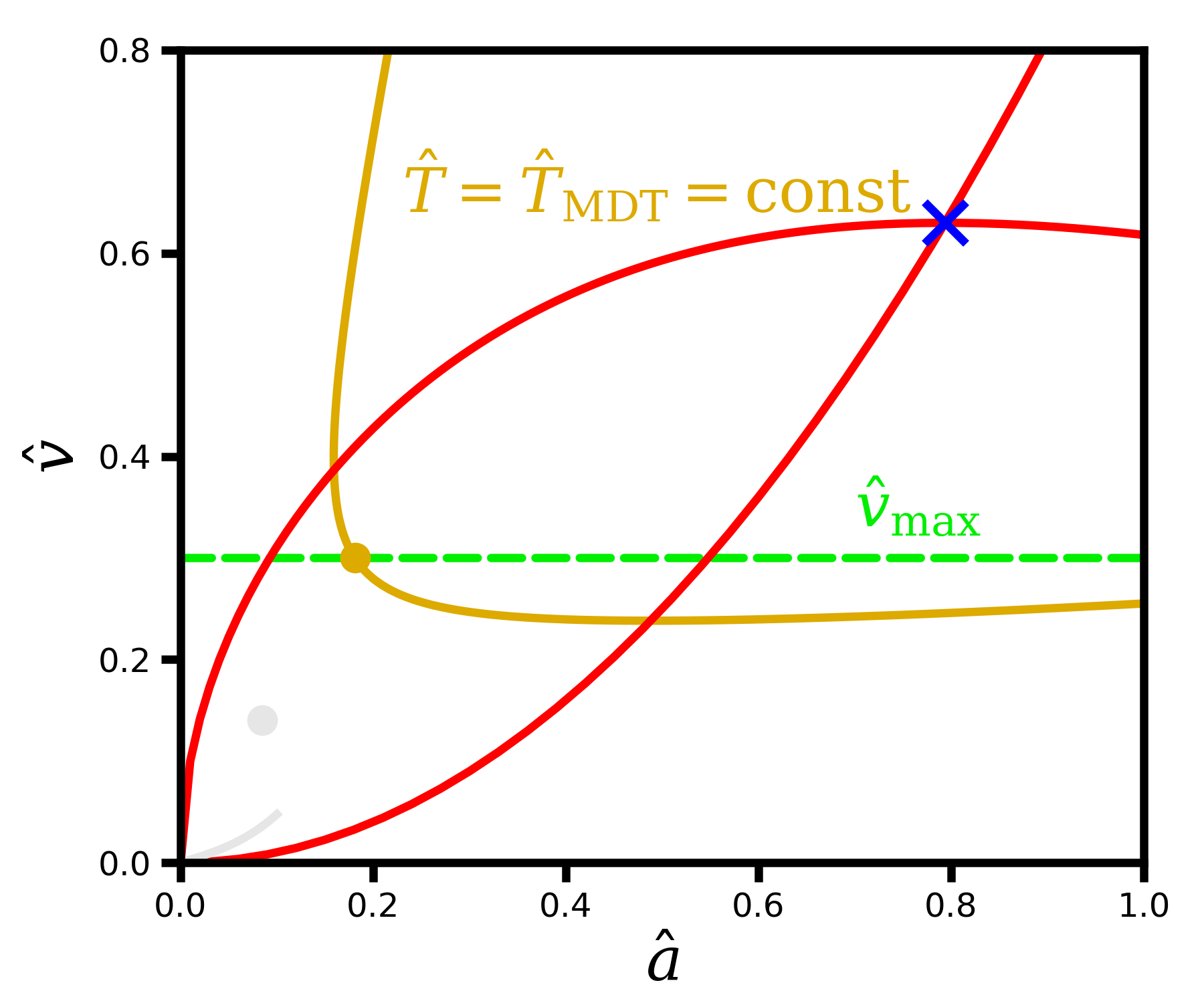}
	\caption{Analogously depiction of MDT with $M=2$ (minimum-acceleration motions) as in Figure~\ref{Figure6}. The left panel shows the case where the solution is determined by the $\hat{a}_{\min}$ curve, whereas the right panel displays the case where the solution is determined by the $\hat{v}_{\max}$ bound.} \label{Figure7}
\end{figure}

Let us now take a closer look at MTT. The four possible cases discussed in Section~\ref{3.2} can be easily recovered and reinterpreted within the diving fish diagram. As above, $\hat{v}_{\max}$ and $\hat{a}_{\max}$ are horizontal and vertical cut-offs, respectively. Then:
\begin{itemize}
	\item 4-phase case: $\hat{v}_{\max} \geq \hat{v}^{\star}$ and $\hat{a}_{\max} \geq \hat{a}^{\star}$ (no active bound), and the solution of MTT is euqivalent to the global minimum (see Figure~\ref{Figure5}, upper left panel). 
	\item 5-phase case: $\hat{v}_{\max} < \hat{v}^{\star}$ and $\hat{a}_{\max} \geq \hat{a}^{\star}$ (velocity bound active), and the solution of MTT is given by the intersection of $\hat{v}_{\max}$ and the $\hat{a}_{\min}$ curve\footnote{Note that the inactive acceleration bound is equivalent to Equation~(\ref{amin}).} (see Figure~\ref{Figure5}, upper right panel).
	\item 6-phase case:  $\hat{v}_{\max} \geq \hat{v}^{\star}$ and $\hat{a}_{\max} < \hat{v}^{\star}$ (acceleration bound active), and the solution of MTT is given by the intersection of $\hat{a}_{\max}$ and the $\hat{v}_{\min}$ curve\footnote{The inactive velocity bound is equivalent to Equation~(\ref{vmin}).} (see Figure~\ref{Figure5}, lower left panel). 
	\item 7-phase case: $\hat{v}_{\max} < \hat{v}^{\star}$ and $\hat{a}_{\max} < \hat{a}^{\star}$ (both bounds active); then $\hat{v} = \hat{v}_{\max}$ and $\hat{a} = \hat{a}_{\max}$ (see Figure~\ref{Figure5}, lower right panel) and Equation~(\ref{5.2T0}) becomes
	\begin{equation}
		T = \frac{s}{v_{\max}} + \frac{v_{\max}}{a_{\max}} + \frac{a_{\max}}{j_{\max}}.
	\end{equation}
\end{itemize}
A slightly different perspective is exhibited in Figure~\ref{Figure4B}, where $\hat{T}$ is displayed as a function of $\hat{T}_1$ and $\hat{T}_2$ (see Equations~(\ref{hatT1}, \ref{hatT2})). The $\hat{v}_{\min}$ curve has become a straight line. While $\hat{T}_2$ is always bounded, $\hat{T}_1$ can become arbitrary high. There is a similar procedure for MDT. Since $\hat{T}_{\mathrm{MDT}}$ is given we work on a fixed contour curve.
\begin{itemize}
	\item minimum-velocity motion $(M=1)$: There are two cases, namely:
	\begin{enumerate}[label=\Alph*)] 
		\item The solution is given by the intersection of the fixed contour curve with the $\hat{v}_{\min}$ curve. This is the case if $\hat{a}_{\max}$ is not active (see Figure~\ref{Figure6}, left panel).
		\item Otherwise, the solution is determined by the intersection of the fixed contour curve and the vertical cut-off defined by $\hat{a}_{\max}$ (see Figure~\ref{Figure6}, right panel).
	\end{enumerate}
	\item minimum-acceleration motion $(M=2)$: The velocity and the acceleration interchanges roles, and there are two cases as well:
	\begin{enumerate}[label=\Alph*)] 
		\item The solution is given by the intersection of the fixed contour curve with the $\hat{a}_{\min}$ curve. This is the case if $\hat{v}_{\max}$ is inactive (see Figure~\ref{Figure7}, left panel).
		\item Otherwise, the solution is determined by the intersection of the fixed contour curve and the horizontal cut-off defined by $\hat{v}_{\max}$ (see Figure~\ref{Figure7}, right panel).
	\end{enumerate}
\end{itemize}
All in all, the relevant optimal values can be obtained by targeted solving of algebraic equations. In addition, the diving fish diagram provides an indication of which case actually occurs, i.e. of how many phases the motion consists and which particular bounds are active. Nevertheless, independent of the cases and MTT or MDT, Equation~(\ref{5.2T0}) remains overall valid.

After the principal ideas are shown for the example of a jerk-controlled motion, we are able to extend those to the next order.
\subsection{Snap-Controlled Motions}\label{5.3}
Following the terminology of~\cite{Biagiotti}, the derivative of the jerk is called snap. We introduce $p_{\max} = w_4$ and $p = w_4 = \max\limits_{t\in[0,T]}|x^{(4)}(t)|$. Equation~(\ref{Eq1Def4}) of Theorem~\ref{theo5} reads now
\begin{align}
	T_0 &= T = \frac{s}{v} + \frac{v}{a} + \frac{a}{j} + \frac{j}{p} \label{5.3T0}\\
	T_1 &= \frac{v}{a} +\frac{a}{j} + \frac{j}{p} \\
	T_2 &= \frac{a}{j} + \frac{j}{p} \\
	T_3 &= \frac{j}{p}
\end{align} 

\begin{figure}[h!]
	\centering
	\includegraphics[width=0.5\textwidth]{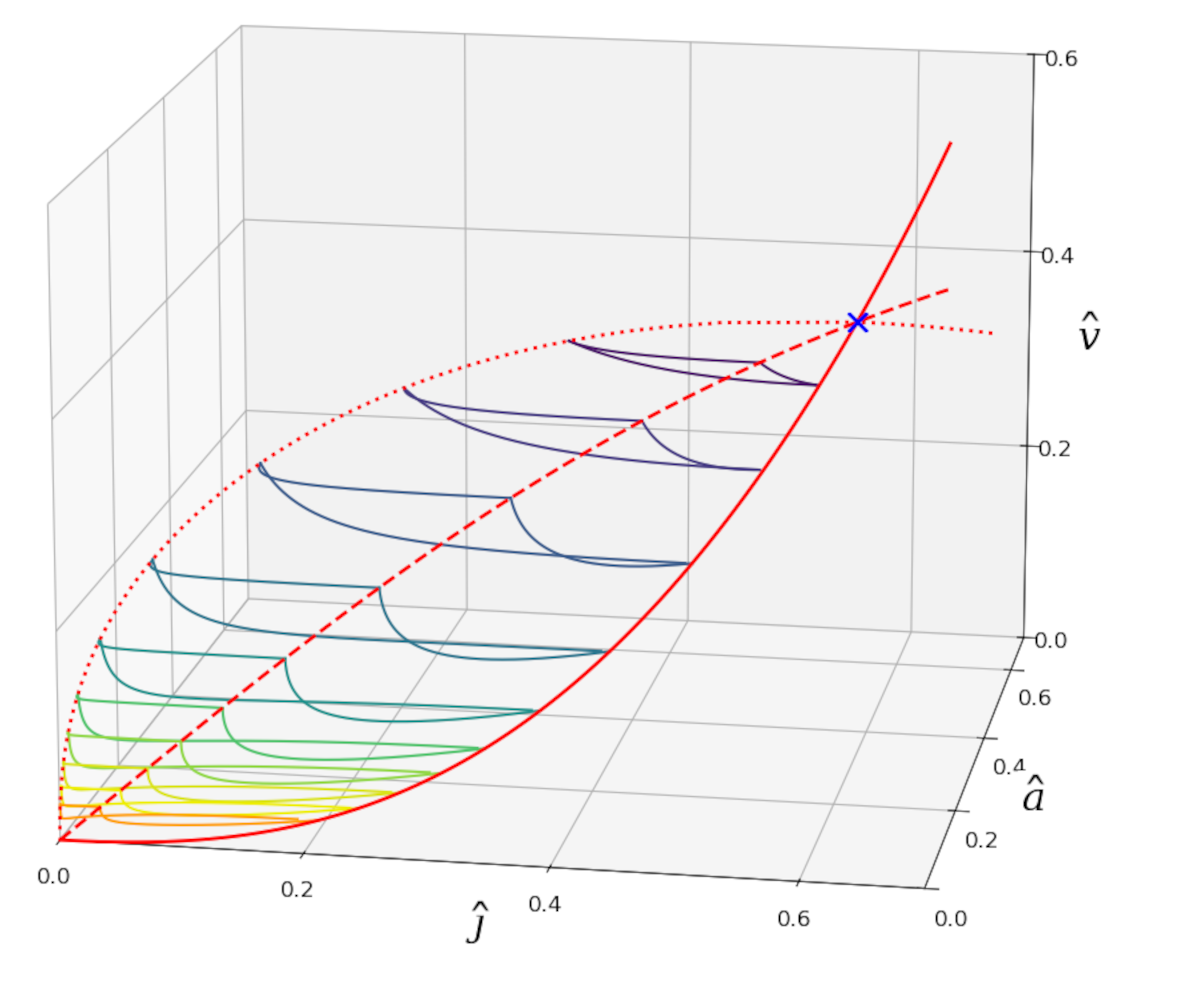}
	\caption{Three-dimensional plot of $\hat{T}$ as a function of $\hat{v}$, $\hat{a}$ and $\hat{\jmath}$ according to Equations~(\ref{5.3T0},~\ref{5.3hatT}). The surfaces of constant $\hat{T}$ value are indicated by the colored curves (from orange for high $\hat{T}$ to deep purple for low $\hat{T}$), where the actual surface embraced by the curves is not shown for a better visibility of the entire figure. The three red curves frame the region of all admissible combinations of ($\hat{v}, \hat{a}, \hat{\jmath}$). They are obtained by combining two of the three conditions of~(\ref{5.3cond}), e.g., combining $\hat{T}_1 = 2\hat{T}_2$ and $\hat{T}_2 = 2\hat{T}_3$ fixes the red solid curve. A single condition like $\hat{T}_1 = 2\hat{T}_2$ corresponds to a surrounding surface (not shown). The global minimum which is determined by the simultaneous consideration of all conditions of~(\ref{5.3cond}) is marked by the blue cross.} \label{Figure8}
\end{figure}

The dimensional transformation
\begin{equation}
	\begin{aligned}
	\hat{T} &=Tp_{\max}^{\frac{1}{4}}s^{-\frac{1}{4}} \\
	\hat{v} &=v p_{\max}^{-\frac{1}{4}}s^{-\frac{3}{4}}\\
	\hat{a} &=ap_{\max}^{-\frac{1}{2}}s^{-\frac{1}{2}} \\ 
	\hat{\jmath} &=j p_{\max}^{-\frac{3}{4}}s^{-\frac{1}{4}} \\
	\hat{s} &= \hat{p} =1 
	\end{aligned} \label{5.3hatT}
\end{equation}
allows us to plot a three-dimensional diving fish diagram (see Figure~\ref{Figure8}). The contour surfaces of constant $\hat{T}$ values are indicated by the colored curves circulating around the fish. The conditions
\begin{equation}
	\hat{T}_0 = 2\hat{T}_1 \quad \hat{T}_1 = 2\hat{T}_2 \quad \hat{T}_2
 = 2\hat{T}_3 \label{5.3cond}
 \end{equation}
 define the three boundary surfaces. Combining two of these conditions yields the red boundary curves. And the combination of all the three conditions fixes the global minimum (see Lemma~\ref{lem9}) located at
 \begin{equation}
 	\hat{\jmath}^{\star} = \sqrt[4]{\frac{1}{8}}, \quad 
 	\hat{a}^{\star} = \sqrt{\frac{1}{8}}, \quad 
 	\hat{v}^{\star} = \sqrt[4]{\frac{1}{32}}, \quad 
 	\hat{T}^{\star} = \sqrt[4]{512}.
 \end{equation}
 We will refrain here from discussing each single option (see~\cite{Biagiotti} for a case-by-case list) and rather incomprehensible drawings of complex intersections but it is easy to get an overview:
 \begin{itemize}
 	\item The bounds $\hat{v}_{\max}$, $\hat{a}_{\max}$ and $\hat{\jmath}_{\max}$ are now cutting planes.
 	\item If no bound is active (i.e. the solution of MTT is equivalent to the global minimum), $\dddot{x}(t)$ is ``fully triangular'', and none of the orders has a cruising phase. Hence this kind of motion consists of 8 phases which is the minimum number of phases for a snap-controlled motion.
 	\item The opposite case -- all bounds are active -- consists of 15 phases. Therefore, there are 8 cases to be distinguished.
 	\item This distinction can be made systematically by switching the conditions of~(\ref{5.3cond}) on or off independently from each other and counting how many combinations result (indeed, $2^3 = 8$).
 \end{itemize}
 If one incorporates yet another order into analysis, i.e. the click $x^{(5)}$ as the derivative of the snap, graphical representation reaches its limits since unfortunately the fourth dimension and hyper-fishes defies appropriate graphical representations.
\subsection{Systematics for Arbitrary Order}\label{5.4}
Based on Theorem~\ref{theo5},
\begin{equation}
	T_n = \sum_{k=n}^{N-1} \frac{x_k}{x_{k+1}}
\end{equation}
for $0\leq n \leq N-1$, the conditions
\begin{equation}
	T_n = 2T_{n+1} \label{5.42T}
\end{equation}
for $0\leq n \leq N-2$ define the hyper-surfaces framing the region of all possible motions. A motion of order $N$ consists at least of $2^{N-1}$ phases. Equation~(\ref{5.42T}) includes $N-1$ conditions with $2^{N-1}$ combinations for switching on and off. This leads to at most $2^{N}-1$ phases for an $N$-th order motion.

We collect some specific examples. Starting with the jerk-controlled motion and
\begin{equation}
	 T = \frac{s}{v} + \frac{v}{a} + \frac{a}{j}
\end{equation}
with at least 4 phases and at most 7 phases, the next order of snap-controlled motions obeys
\begin{equation}
	T = \frac{s}{v} + \frac{v}{a} + \frac{a}{j} + \frac{j}{p}.
\end{equation}
As already mentioned, such a motion has at least 8 phases and at most 15 phases. The click already mentioned at the end of Section~\ref{5.3} is the 5th order $x^{(5)}(t)$. Setting $c = x_5 = \max\limits_{t\in[0,T]}|x^{(5)}(t)|$ as achieved click leads to
\begin{equation}
		T = \frac{s}{v} + \frac{v}{a} + \frac{a}{j} + \frac{j}{p} + \frac{p}{c}.
\end{equation}
Also even higher orders, for example the yeet $x^{(6)}(t)$ and the pow $x^{(7)}(t)$, can be considered (\cite{Boryga} goes to 9th order). We will continue with an arbitrary order $N$.
The dimensional transformation to obtain $\hat{x}_0 = \hat{x}_N = 1$ is in general given by
\begin{align}
	\hat{x}_n &= x_n x_0^{\frac{n}{N}-1} x_N^{-\frac{n}{N}} \label{xxx}\\
	\hat{T} &= T x_0^{-\frac{1}{N}} x_N^{\frac{1}{N}}.
\end{align}
 The global minimum is determined by considering Equation~(\ref{5.42T}) as a system of $N-1$ equations.
 \begin{lemma}\label{lem9}
 	The global minimum of an $N$-th order motion is located at 
 	\begin{align}
 		\hat{x}^{\star}_{N-1} &= 2^{\frac{1}{2N}(N-1)(2-N)} \label{L9.1}\\
 		\hat{x}^{\star}_n &= \hat{x}_{N-1}^{{\star}{N-n}} 2^{\frac{1}{2}(n-N+1)(n-N+2)} = 2^{\frac{n}{2}(\frac{2}{N}-N+n)} \label{L9.2}
 \end{align}
 with the minimal time horizon
 \begin{equation}
 	\hat{T}^{\star}_0 = 2^{\frac{1}{2N}(N-1)(N+2)}. \label{L9.3}
 \end{equation}
 In addition, one has
 \begin{equation}
 	\hat{T}^{\star}_n = 2^{N-n-1} \hat{x}^{\star}_{N-1}. \label{L9.4}
 \end{equation}
 \end{lemma}
 \begin{proof}
 	Clearly, $\hat{T}^{\star}_{n-1} = \hat{x}^{\star}_{N-1}$. Due to Equation~(\ref{5.42T}) it is easy to find that
 	\begin{equation}
 		\hat{T}^{\star}_n = 2^{N-n-1} \hat{x}^{\star}_{N-1}
 	\end{equation}
 	valid for $0 \leq n \leq N-1$. Combining Equation~(\ref{L9.4}) with Equation~(\ref{Arrive}) (see Lemma~\ref{Lem8} as well) we arrive at
 	\begin{equation}
 		\hat{T}^{\star}_{N-1} - \hat{T}^{\star}_n = \hat{x}^{\star}_{N-1} 2^{N-n-1}
 	\end{equation}
 	and 
 	\begin{equation}
 		\hat{x}^{\star}_{n-1} = \hat{x}^{\star}_{N-1} \hat{x}^{\star}_n 2^{N-n-1}.
 	\end{equation}
 	With the ansatz
 	\begin{equation}
 		\hat{x}^{\star}_n = \hat{x}_{N-1}^{{\star}{N-n}} 2^{h(n)}
 	\end{equation}
 	and a function $h$, one obtains
 	\begin{equation}
 		h(n-1) - h(n) = N-n-1
 	\end{equation}
 as a recursion for $h$ (with $h(N-1)=0$). Then the solution is $h(n) = \frac{1}{2}(n-N+1)(n-N+2)$ which proves Equation~(\ref{L9.2}). The condition $\hat{T}^{\star}_0 = 2\hat{T}^{\star}_1$ which is equivalent to $\frac{1}{\hat{x}^{\star}_1} = \hat{T}^{\star}_1$ finally determines $\hat{x}^{\star}_{N-1}$ by using Equations~(\ref{L9.2}) and (\ref{L9.4}). Equation~(\ref{L9.4}) is a direct implication.
 \end{proof}
 The Table~\ref{Tab1} lists some examples.
 \begin{table}
 	 	\centering
  	 \begin{tabular}{ccccc} 
 	N & $\hat{T}^{\star}_0$ & $\hat{v}^{\star}$ & $\hat{a}^{\star}$  & $\hat{\jmath}^{\star}$ \\ 
  	\hline  \vspace{3pt}
 	2 & 2 & 1 &  &  \\
 		\vspace{3pt}
 	3 & $\sqrt[3]{32}$ & $\sqrt[3]{\frac{1}{4}}$ & $\sqrt[3]{\frac{1}{2}}$ &  \\
 		\vspace{3pt}
 	4 & $\sqrt[4]{512}$ & $\sqrt[4]{\frac{1}{32}}$ & $\sqrt{\frac{1}{8}}$ & $\sqrt[4]{\frac{1}{8}}$ \\ 
  	\hline
 \end{tabular}
 	\caption{Loci of the global minimum for the orders $N=2,3,4$.}\label{Tab1}
  \end{table}
 It can be seen from Equation~(\ref{L9.3}) that the minimal time grows exponentially with the order of the motion\footnote{While for the jerk-controlled motion $\hat{T}^{\star}_0 \approx 3.17$, the case $N=6$ (yeet-controlled motion) and $N=7$ (pow-controlled motion) yield  $\hat{T}^{\star}_0 \approx 10.08$ and  $\hat{T}^{\star}_0 \approx 14.49$.}. So the time intervals, where the highest order is non-vanishing, become rapidly shorter and the ratios $\frac{x_n}{x_{n+1}}$ for higher $n$ shrink.
 
\subsection{Computation of Minimum-Time Trajectories}\label{5.5}
Let us eventually consider MTT algorithmically. Although not crucial, we will mainly work with dimensionless quantities. It is convenient to define the earliest time point of achievement $T_n$ (see Definition~\ref{Def4}) as a function on $\mathbb{R}^{N+1}$: For $z=(z_0,\ldots,z_N)\in \mathbb{R}^{N+1}$ set
\begin{equation}
	T_n(z) = \sum_{k=n}^{N-1} \frac{z_k}{z_{k+1}}. \label{5.5E1}
\end{equation}
The main idea\footnote{Note that minimizing the time horizon is a convex optimization problem since the objective function $T_0$ is convex and the convex combination of two admissible trajectories is admissible as well (cf.~\cite{boyd} for instance).} is based on the observation that the function $T_0(z)$ attains no inner minimum\footnote{Regarding the entire plane $\mathbb{R}^{N+1}$, it is a straight forward calculation to show that $z_0 = \ldots = z_N = 1$ is the global minimum.} within the admissible region framed by the hypersurfaces encoded in $T_{n-1}(z) = 2T_n(z)$. Therefore, either a bound $w_n$ is active, or the equation $T_{n-1}(z) = 2T_n(z)$ holds.

\begin{figure}[h!]
	\centering
	\includegraphics[width=0.5\textwidth]{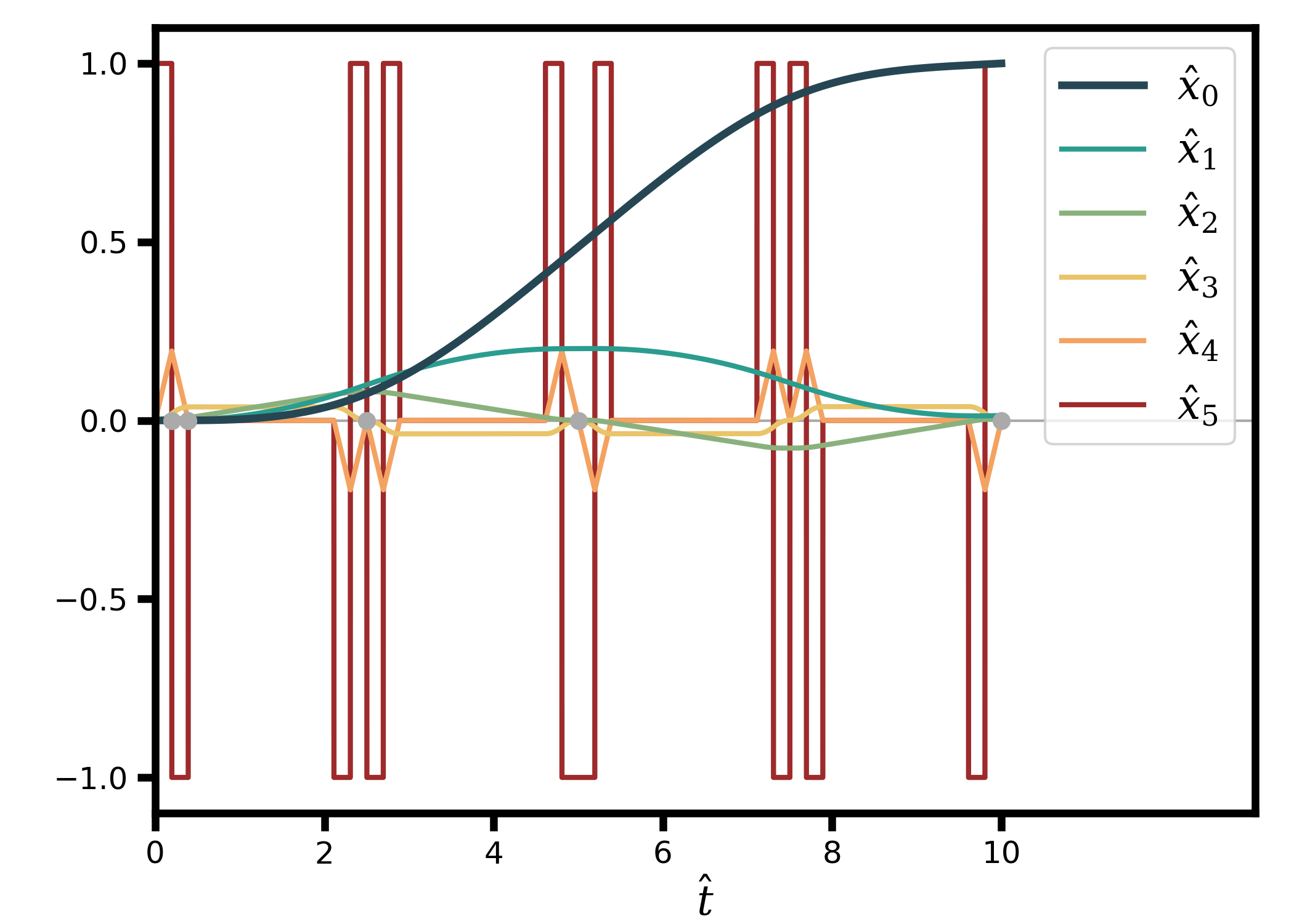}
	\caption{Dimensionless profiles of the minimum-time trajectory for $s= 20$~m (see Table~\ref{Tab2}), where $\hat{x}_0$ show the position as a function of the time (higher orders accordingly). The dimensionless time horizon of the motion is $\hat{T}_0 = 10.00$ (equivalent to 11.49~s). The maximal velocity is initially achieved at $\hat{T}_1 = 5.00$ (5.74~s) -- and higher orders analogously.} \label{Figure9}
\end{figure}

\begin{figure}[h!]
	\centering
	\includegraphics[width=0.5\textwidth]{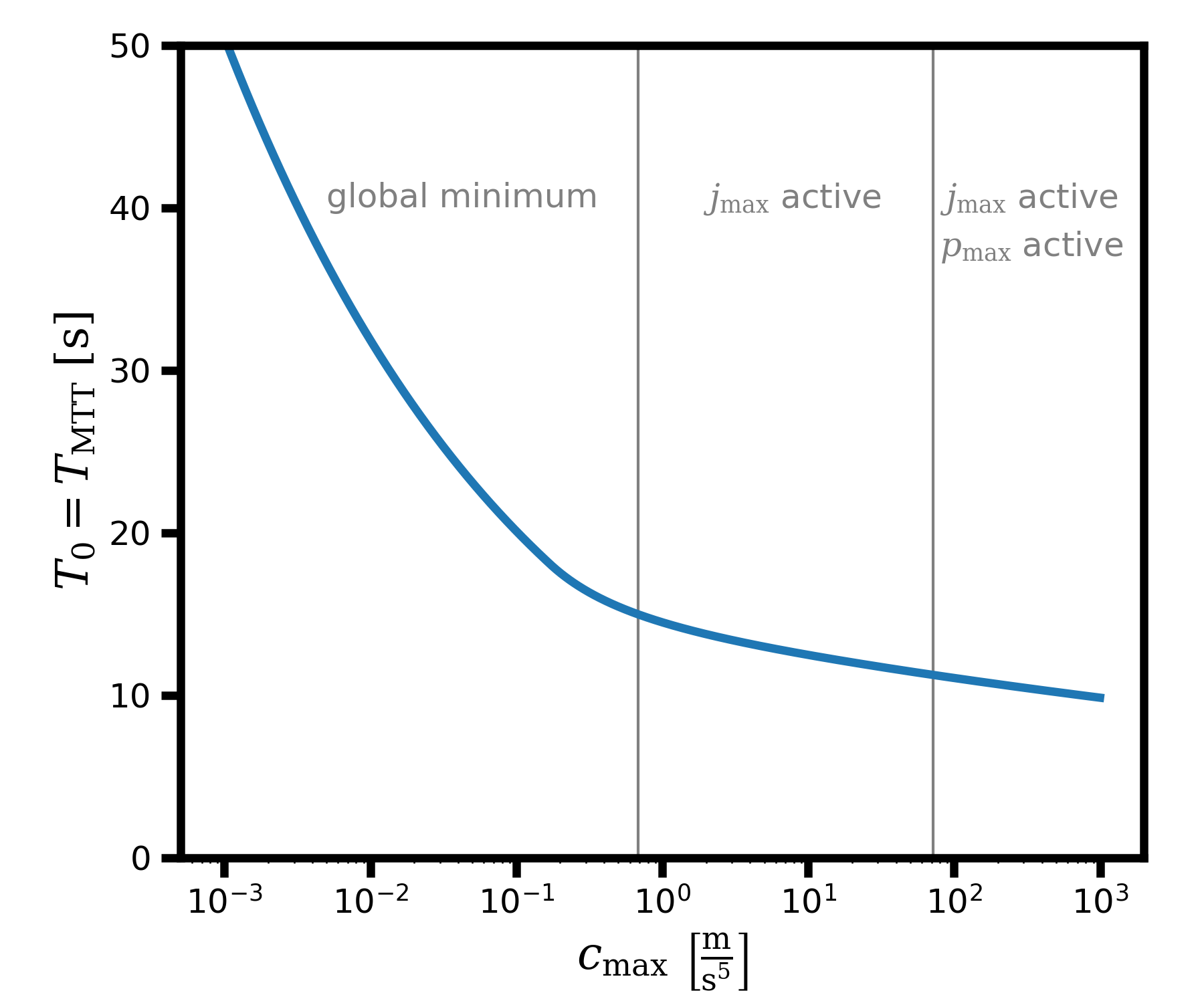}
	\caption{Time horizon of the minimum-time motion as a function of $c_{\max}$. The other parameters are of the example displayed in Figure~\ref{Figure9}.} \label{Figure10}
\end{figure}

Let a distance $s$ as well as the bounds $w_1,\ldots, w_N$ be given. Then the steps of the algorithm for finding the minimum-time motion are as follows:
\begin{enumerate}
	\item Transform $w_n$ to $\hat{w}_n$ according to Equation~(\ref{xxx}), and set $\hat{u} = (\hat{u}_0, \ldots, \hat{u}_N)$ with
	\begin{equation}
		\hat{u}_n = \begin{cases}
			1 &n=0,N\\
			\min\left\lbrace \hat{w}_n, \hat{x}_n^{\star}\right\rbrace  &0<n<N
		\end{cases}
	\end{equation}
	\item Collect all indices $1 \leq n \leq N-1$ with $T_{n-1}(\hat{u}_n) < 2T_n(\hat{u}_n)$ (see Equation~(\ref{5.5E1})), i.e. find the set 
	\begin{equation}
		B = \left\lbrace n~|~T_{n-1}(\hat{u}_n) < 2T_n(\hat{u}_n) \wedge 1 \leq n \leq N-1\right\rbrace.
	\end{equation}
	\item Solve\footnote{E.g., with Newton's method.} the system of equations
	\begin{equation}
		T_{n-1}(\hat{z}) = 2T_n (\hat z) \quad (n \in B)
	\end{equation}
	for $\hat{z} = (\hat{z}_0, \ldots, \hat{z}_N), \hat{z}_0 = \hat{z}_N = 1$ including $\hat{z}_n = \hat{u}_n$ for $n\notin B$.
	\item Perform a back transformation from $\hat{z}$ to gain the solution $x_{\mathrm{MTT}}$ and the time points of earliest achievement as $T_n(x_{\mathrm{MTT}})$. 
\end{enumerate}
Let us illustrate the algorithm by the following example of 5th order. We seek for the minimum time motion accomplishing $s=20$~m subject to constraints
\begin{align*}
	w_1 &= v_{\max} = 7~\frac{\mathrm{m}}{\mathrm{s}} \\
	w_2 &= a_{\max} = 2~\frac{\mathrm{m}}{\mathrm{s}^2} \\
	w_3 &= j_{\max} = 0.5~\frac{\mathrm{m}}{\mathrm{s}^3} \\
	w_4 &= p_{\max} = 6~\frac{\mathrm{m}}{\mathrm{s}^4} \\
	w_5 &=c_{\max} = 10~\frac{\mathrm{m}}{\mathrm{s}^5}.
\end{align*}
Table~\ref{Tab2} summarizes the steps to determine the minimum-time motion. The ratio $\frac{T_2(\hat{u})}{T_3(\hat{u})}>2$ fixes $\hat{z}_3=\hat{w}_3=0.0379$ (and $j_{\mathrm{MTT}} = j_{\max} = 0.5~\frac{\mathrm{m}}{\mathrm{s}^3}$). For the remaining quantities $\hat{z}_n$ we have to solve
\begin{align}
	\frac{1}{\hat{z}_1} &= \frac{\hat{z}_1}{\hat{z}_2} + \frac{\hat{z}_2}{\hat{w}_3} + \frac{\hat{w}_3}{\hat{z}_4} + \hat{z}_4 \label{5.5E2}\\
	\frac{\hat{z}_1}{\hat{z}_2} &= \frac{\hat{z}_2}{\hat{w}_3} + \frac{\hat{w}_3}{\hat{z}_4} + \hat{z}_4 \label{5.5E3}\\
	\frac{\hat{w}_3}{\hat{z}_4} &= \hat{z}_4 \label{5.5E4}
\end{align}
\begin{table}
	\centering
	\begin{tabular}{K{1cm}K{0.9cm}K{0.9cm}K{0.9cm}K{0.9cm}K{0.9cm}K{0.9cm}K{0.9cm}K{0.9cm}K{0.9cm}K{0.9cm}K{0.9cm}K{0.9cm}}
	\hline
	& \mc{$s$} & \mc{$v_{\max}$} & \mc{$a_{\max}$}  & \mc{$j_{\max}$} & \mc{$p_{\max}$} & \mc{$c_{\max}$} \\ 
	&&& \mc{$w_1$} & \mc{$w_2$} & \mc{$w_3$} & \mc{$w_4$} & \mc{$w_5$} \\
	&\mc{20~m} &\mc{$7~\frac{\mathrm{m}}{\mathrm{s}}$} &\mc{$2~\frac{\mathrm{m}}{\mathrm{s^2}}$} & \mc{$0.5~\frac{\mathrm{m}}{\mathrm{s^3}}$} & \mc{$6~\frac{\mathrm{m}}{\mathrm{s^4}}$} & \mc{$10~\frac{\mathrm{m}}{\mathrm{s^5}}$} \\
	\hline
	\multirow{2}{*}{step 1} &&&\mc{$\hat{w}_1$} & \mc{$\hat{w}_2$}  & \mc{$\hat{w}_3$} & \mc{$\hat{w}_4$}\\
	&&&\mc{0.4020} & \mc{0.1320}  & \mc{0.0379} & \mc{0.5223}\\
	& \mc{$\hat{u}_0$} & \mc{$\hat{u}_1$} & \mc{$\hat{u}_2$}  & \mc{$\hat{u}_3$} & \mc{$\hat{u}_4$} & \mc{$\hat{u}_5$} \\
	& \mc{1} & \mc{0.2873} & \mc{0.1320} & \mc{0.0379} & \mc{0.4353} & \mc{1} \\
	\hline
	\multirow{4}{*}{step 2} & \mc{$T_0(\hat{u})$} & \mc{$T_1(\hat{u})$} & \mc{$T_2(\hat{u})$}  & \mc{$T_3(\hat{u})$} & \mc{$T_4(\hat{u})$} &&\\
	& \mc{9.6631} & \mc{6.1809} & \mc{4.0045} & \mc{0.5223} & \mc{0.4353} &&  \\
	&&\mc{ratio} &\mc{ratio} &\mc{ratio} &\mc{ratio}&& \\
	&&\mc{1.56}  &\mc{1.54}  &\mc{7.67}  &\mc{1.20}&&\\
	\hline
	\multirow{2}{*}{step 3} & \mc{$\hat{z}_0$} &\mc{$\hat{z}_1$} & \mc{$\hat{z}_2$} & \mc{$\hat{z}_3$} & \mc{$\hat{z}_4$} & \mc{$\hat{z}_5$} \\
	& \mc{1} & \mc{0.1999} & \mc{0.0799} & \mc{0.0379} & \mc{0.1947} & \mc{1} \\
	\hline
	\multirow{4}{*}{step 4}& \mc{$x_{\mathrm{MTT},0}$} & \mc{$x_{\mathrm{MTT},1}$} &	\mc{$x_{\mathrm{MTT},2}$} & \mc{$x_{\mathrm{MTT},3}$} & \mc{$x_{\mathrm{MTT},4}$} & \mc{$x_{\mathrm{MTT},5}$}\\
	& \mc{20~m} & \mc{3.48~$\frac{\mathrm{m}}{\mathrm{s}}$} & \mc{1.21~$\frac{\mathrm{m}}{\mathrm{s^2}}$} & \mc{0.5~$\frac{\mathrm{m}}{\mathrm{s^3}}$} & \mc{2.24~$\frac{\mathrm{m}}{\mathrm{s^4}}$} & \mc{10~$\frac{\mathrm{m}}{\mathrm{s^5}}$}\\
	&\mc{$T_0$}	&\mc{$T_1$} &\mc{$T_2$}	&\mc{$T_3$} &\mc{$T_4$}\\
	&\mc{11.49~s}&\mc{5.74~s}&\mc{2.87~s}&\mc{0.45~s}&\mc{0.22~s}\\
	\hline
\end{tabular}
\caption{Summary of exemplary time minimization for a distance $s=20$~m and the above constraints. The minimal time is $T=8.4481$~s.}\label{Tab2}
\end{table}
where we already inserted $\hat{w}_3$ as $\hat{z}_3$. Starting from~(\ref{5.5E4}), the solution is not difficult to find (see Table~\ref{Tab2}). We have $\hat{T}_0 = 10.0$ corresponding to $T = T_0 = 11.49~\mathrm{s}$ as minimal time. Figure~\ref{Figure9} shows all orders of the optimal trajectory. Note that there are two interesting cases when varying the bound of the highest order $w_5 = c_{\max}$ (see Figure~\ref{Figure10}):
\begin{itemize}
	\item If $c_{\max} = 72~\frac{\mathrm{m}}{\mathrm{s}^5}$, then both $j_{\max}$ and $p_{\max}$ are always reached. Independent of $c_{\max}$, $a_{\max}$ and $v_{\max}$ can never be reached since $s$ is too short.
	\item If $c_{\max} \leq 0.683$, then no bound is reached anymore and the solution of MTT is equivalent to the global minimum $x_1^{\star}, \ldots, x^{\star}_{N-1}$ (see Lemma~\ref{lem9}). Graphically speaking, the small value of $c_{\max}$ catapults the transformed boundaries $\hat{w}_n$ to high values lying for beyond the tail of the fish.
\end{itemize}
In Appendix~\ref{App.A} we provide an implementation of the algorithm in Python and one implementation in C++. Solving the above example with the Python implementation\footnote{Note that we did not place much importance to the best possible high-end implementation. For instance, the Python implementation could be enhanced through a modification of Newton's method by incorporating analytic derivatives instead of using numerical derivatives.} took about 1~ms on an Intel\textsuperscript{\textregistered} Core\textsuperscript{\texttrademark} i5-8265U.

This algorithm allows a direct computation of the reached values of velocity, acceleration, etc. Since the values $T_n$ follow directly the trajectory can be constructed directly as well.

Finally, with the insights gained from the above algorithm, the MDT optimization can be treated similarly. Let us sketch the procedure. Recall that a distance $s$ and a time horizon $T_{\mathrm{MDT}}$ are given. The objective function is $x_M = \max\limits_{t\in[0,T]} \left|x^{(M)} (t) \right| \rightarrow \min$ for $1 \leq M < N$. As a first step, $T_M$ has to be minimized. Observe that the smaller $T_M$ the longer $x_M$ can be held and consequently the smaller becomes the objective function. The minimization of $T_M$ can be done by applying the above algorithm to the reduced problem of the orders $n = N \ldots M$. This yields a part of the solution, namely $x_{\mathrm{MDT}, N}$, \ldots, $x_{\mathrm{MDT}, M}$. It remains to find $x_{\mathrm{MDT}, M-1}$, \ldots, $x_{\mathrm{MDT}, 1}$. To this end, we have to proceed analogously as above but with two differences: first, the part of the solution for $n \geq M$ is already determined by the $T_M$ minimization. Second, $T_0$ is fixed as $T_{\mathrm{MDT}}$.
	\section{Summary and Outlook} \label{Section6}
This paper presents a method for optimal trajectory planning based on one single equation. This equation, for example in the form 
\begin{equation}
	T = \frac{s}{v} + \frac{v}{a} + \frac{a}{j}
\end{equation}
for jerk-controlled motions, is valid for any number of phases the motion can have. Furthermore, arbitrarily high orders can be considered. For example, 
\begin{equation}
	T = \frac{s}{v} + \frac{v}{a} + \frac{a}{j} + \frac{j}{p} + \ldots = \sum_{n=0}^{N-1} \frac{x_n}{x_{n+1}}.
\end{equation}
This allows for a reinterpretation and solution of well-known standard trajectory optimization problems (such as motions with minimum time or motions with minimum velocity). In particular, the condition 
\begin{equation}
	T_n = 2T_{n+1}
\end{equation} 
enables the determination of the active constraints for the respective case. Additionally, the effect of incorporating higher orders on the time horizon can be considered as well.
Therefore, it seems surprising that -- although the equation is not exceptionally complex and covers a large number of cases -- it has not previously been presented in this form in the literature.

The equations for $T_n$ in combination with the conditions $T_n = 2T_{n+1}$ opens up the possibility of direct calculation of the minimum-time motion.

Besides the practical case of application, it should be noted that motions with the minimum velocity are also energy-optimal in many cases, because a long phase without acceleration combined with a short phase of higher acceleration outweights motions with long phases of lower acceleration.

Clearly, many extensions exist. First, the zero boundary condition can be relaxed. Furthermore, it would be possible to consider other classes of trajectories (for example trajectories with circular blends of trigonometric trajectories). The expected results would then likely be considerably more cumbersome than those obtained here.

Another possibility -- particularly mathematically interesting -- lies in the limit $N \rightarrow \infty$ of infinite order. This raises questions of convergence (see the series in Equation~(\ref{EqTheo5}) in Theorem~\ref{theo5}) as well as questions concerning the relationship with $\mathcal{C}^{\infty}$ functions and mollifiers. Due to the infinitely large number of intervals in a finite time interval, it is also worthwhile to consider methods of fractal geometry (such as fractal strings, see~\cite{Lapidus}). However, all these investigations of interesting connections are left for future work.
	\section*{Acknowledgement}
	The work of Rico Zöllner was supported by DFG project no.~430149671. The author thanks Frank Schulze and Ella Jannasch for their support and valuable comments.
\begin{appendices}
\renewcommand{\appendixtocname}{Appendices}
\renewcommand{\appendixpagename}{Appendices}
	\section{Implementations of the Optimization Algorithm} \label{App.A}
Both implementations below take the transformed boundaries as input and return the solution as dimensionless quantities.
\subsection{Python Implementation}
{\footnotesize
\begin{lstlisting}[label={lst:listing-cpp}, language=Python]
import numpy as np
from scipy.optimize import root
 
# ------------------------------------------------------------
# Step 1: computation of vector u
# ------------------------------------------------------------
def compute_u(N, w):
    w = np.asarray(w, dtype=float)
 
    u = np.zeros(N + 1, dtype=float)
    u[0] = 1.0
    u[N] = 1.0
 
    for n in range(1, N):
        exponent = 0.5 * n * (2.0 / N - N + n)
        bound = 2 ** exponent
        u[n] = min(w[n - 1], bound)
 
    return u
 
# ------------------------------------------------------------
# Def. of T_n: T_n(x) = Sum_{k=n}^{N-1} x_k/x_{k+1}
# ------------------------------------------------------------
def T(n, x):
    N = len(x) - 1
    return np.sum(x[n:N] / x[n+1:N+1])
 
# ------------------------------------------------------------
# Step 2: determination of the set B
# ------------------------------------------------------------
def compute_B(z):
    N = len(z) - 1
    B = []
 
    for n in range(1, N): 
        if T(n - 1, z) < 2 * T(n, z):
            B.append(n)
 
    return np.array(B, dtype=int)
 
# ------------------------------------------------------------
# Step 3: setting up the equation system
# ------------------------------------------------------------
def residuals(x, z, B):
    N = len(x) - 1
    eqs = []
 
    B = set(B)
 
    # x_0 = x_N = 1
    eqs.append(x[0] - 1)
    eqs.append(x[N] - 1)
 
    # equations for n=1, ..., N-1
    for n in range(1, N):
        if n in B:
            eqs.append(T(n - 1, x) - 2*T(n, x))  # n in B
        else:
            eqs.append(x[n] - z[n])  # n not in B
 
    return np.array(eqs)
 
# ------------------------------------------------------------
# main algorithm
# ------------------------------------------------------------
def solve_system(N, w, x0=None):
 
    u = compute_u(N, w)
    B = compute_B(u)
 
    if x0 is None:
        x0 = u.copy()
 
    sol = root(residuals, x0, args=(u, B), method='hybr')
 
    if not sol.success:
        raise RuntimeError("no solution", sol.message)
 
    x = sol.x
 
    if np.any(x <= 0):
        raise RuntimeError("solution not positive")
 
    t = np.array([T(n, x) for n in range(0, N)])
 
    return x, t, B, u
 
# ------------------------------------------------------------
# example
# ------------------------------------------------------------
if __name__ == "__main__":
    N = 5
    w = np.array([0.402044,0.131951,0.037893,0.522330])
 
    x, t, B, u = solve_system(N, w)
 
    print("u =", u)
    print("B =", B)
    print("x =", x)
    print("t =", t)
\end{lstlisting}
}
\subsection{C++ Implementation}
{\footnotesize
\begin{lstlisting}[label={lst:listing-cpp}, language=C++]
#include <iostream>
#include <vector>
#include <cmath>
#include <stdexcept>
#include <iomanip>
 
using namespace std;
 
// ------------------------------------------------------------
// Definition of the functions T_n: T_n(x) = Sum_{k=n}^{N-1}  x_k/x_{k+1}
// ------------------------------------------------------------
double T(int n, const vector<double>& x) {
    int N = x.size() - 1;
    double s = 0.0;
    for (int k = n; k <= N - 1; k++) {
        s += x[k] / x[k + 1];
    }
    return s;
}
 
// ------------------------------------------------------------
// Step 1: computation of vector u
// ------------------------------------------------------------
vector<double> compute_u(int N, const vector<double>& w) {
 
    vector<double> u(N + 1);
    u[0] = 1.0;
    u[N] = 1.0;
 
    for (int n = 1; n <= N - 1; n++) {
        double exponent = 0.5 * n * (2.0 / N - N + n);
        double bound = pow(2.0, exponent);
        u[n] = min(w[n - 1], bound);
    }
    return u;
}
 
// ------------------------------------------------------------
// Step 2: determination of the set B
// ------------------------------------------------------------
vector<int> compute_B(const vector<double>& z) {
    int N = z.size() - 1;
    vector<int> B;
 
    for (int n = 1; n <= N - 1; n++) {
        if (T(n - 1, z) < 2.0 * T(n, z))
            B.push_back(n);
    }
    return B;
}
 
// ------------------------------------------------------------
// Newton method: computation of the residuals
// ------------------------------------------------------------
vector<double> residuals(const vector<double>& x,
                         const vector<double>& z,
                         const vector<int>& Bset) {
 
    int N = x.size() - 1;
    vector<double> r(N + 1, 0.0);
 
    // B as boolean index
    vector<bool> inB(N + 1, false);
    for (int b : Bset) inB[b] = true;
 
    // x_0 = x_N = 1
    r[0] = x[0] - 1.0;
    r[N] = x[N] - 1.0;
 
    // equations
    for (int n = 1; n <= N - 1; n++) {
        if (inB[n]) {  // n in B
            r[n] = T(n - 1, x) - 2.0 * T(n, x);
        } else {  // n not in B
            r[n] = x[n] - z[n];
        }
    }
    return r;
}
 
// ------------------------------------------------------------
// Newton method for solving the equation system
// ------------------------------------------------------------
vector<double> solve_nonlin(const vector<double>& z, const vector<int>& B) {
    int N = z.size() - 1;
    vector<double> x = z;        // initial value: u
    const double eps = 1e-10;
    const double h = 1e-6;
    const int maxIter = 200;
 
    for (int it = 0; it < maxIter; it++) {
 
        vector<double> F = residuals(x, z, B);
 
        // checking for convergence
        double norm = 0.0;
        for (double v : F) norm += v * v;
        norm = sqrt(norm);
 
        if (norm < eps) {
            return x;
        }
 
        // numerical Jacobian
        vector<vector<double>> J(N + 1, vector<double>(N + 1, 0.0));
 
        for (int j = 0; j <= N; j++) {
            vector<double> xh = x;
            xh[j] += h;
            vector<double> Fh = residuals(xh, z, B);
 
            for (int i = 0; i <= N; i++) {
                J[i][j] = (Fh[i] - F[i]) / h;
            }
        }
 
        // solving linear system J * dx = -F with Gauss
        vector<double> dx = F;
        for (double &v : dx) v = -v;
 
        // Gauss elimination:
        for (int i = 0; i <= N; i++) {
 
            double piv = J[i][i];
            if (fabs(piv) < 1e-14)
                throw runtime_error("Newton: singular matrix.");
 
            // normalizing
            for (int j = i; j <= N; j++) J[i][j] /= piv;
            dx[i] /= piv;
 
            // elimination
            for (int k = 0; k <= N; k++) {
                if (k == i) continue;
                double f = J[k][i];
                for (int j = i; j <= N; j++) {
                    J[k][j] -= f * J[i][j];
                }
                dx[k] -= f * dx[i];
            }
        }
 
        // next iteration point
        for (int i = 0; i <= N; i++) {
            x[i] += dx[i];
        }
    }
 
    throw runtime_error("Newton method did not converge.");
}
 
// ------------------------------------------------------------
// main algorithm and example
// ------------------------------------------------------------
int main() {
 
    int N = 5;
    vector<double> w = {0.402044,0.131951,0.037893,0.522330};
 
    vector<double> u = compute_u(N, w);
    vector<int> B = compute_B(z);
 
    vector<double> x = solve_nonlin(u, B);
 
    // t vektor: t_n = T_n(x)
    vector<double> t(N);
    for (int n = 0; n < N; n++) {
        t[n] = T(n, x);
    }
 
    // output
    cout << "u = ";
    for (double v : u) cout << v << " ";
    cout << "\nB = ";
    for (int b : B) cout << b << " ";
    cout << "\nx = ";
    for (double v : x) cout << v << " ";
    cout << "\nt = ";
    for (double v : t) cout << v << " ";
    cout << endl;
 
    return 0;
}
\end{lstlisting}
}
\end{appendices}
	\bibliography{Riggns_Cut_Literature}
	\bibliographystyle{tip2.bst}
\end{document}